\documentclass{amsart}

\usepackage{bbm}
  \usepackage{bm}
\usepackage[lite]{amsrefs}
\usepackage{amssymb}
\usepackage[all,cmtip]{xy}
\usepackage{amsmath}
\usepackage{amsthm}
\usepackage{amsfonts}
  \usepackage{thmtools}
\usepackage[hypertexnames=false]{hyperref}
\usepackage{enumerate}
\usepackage[titletoc]{appendix}
\usepackage{fancyhdr}
\usepackage[dvipsnames]{xcolor}
\usepackage{comment}
\usepackage{tikz-cd}
\usetikzlibrary{babel}
\usepackage{tikzpagenodes}
\usepackage{caption}
\usetikzlibrary[patterns]
\usepackage[margin = 1in]{geometry}
\usepackage{setspace}
\usepackage[T1]{fontenc}
  \usepackage{cleveref}
  \crefalias{subsection}{section}

\allowdisplaybreaks

\declaretheorem[name=Theorem,Refname={Theorem,Theorems},parent=section]{thm}

\declaretheorem[name=Lemma,sibling=thm]{lem}
\declaretheorem[name=Proposition,Refname={Proposition,Propositions},sibling=thm]{prop}
\declaretheorem[name=Definition,Refname={Definition,Definitions},sibling=thm]{defn}
\declaretheorem[name=Example,sibling=thm]{ex}
\declaretheorem[name=Remark,sibling=thm]{rem}

\declaretheorem[name=Definition-Theorem,sibling=thm]{defn-thm}
\declaretheorem[name=Notation,sibling=thm]{notation}

   \renewcommand\[{\begin{equation}}
\renewcommand\]{\end{equation}}
\numberwithin{equation}{section}

\usepackage{amsxtra}
\usepackage{graphicx}
\usepackage{bbm,newlfont}
\usepackage{amscd}
\usepackage[english]{babel}
\usepackage{cancel}
\usepackage{amsmath,amsthm,amssymb}
\usepackage{enumerate}
\usepackage{color}
\usepackage{varwidth}
\usepackage{mathrsfs}  
\usepackage{setspace}
\usepackage{collectbox}
\usepackage[utf8]{inputenc} % usually not needed (loaded by default)
\usepackage[T1]{fontenc}


\usepackage{graphicx}% just for the example text

\allowdisplaybreaks

\makeatletter
\tikzset{
  column sep/.code=\def\pgfmatrixcolumnsep{\pgf@matrix@xscale*(#1)},
  row sep/.code   =\def\pgfmatrixrowsep{\pgf@matrix@yscale*(#1)},
  matrix xscale/.code=%
    \pgfmathsetmacro\pgf@matrix@xscale{\pgf@matrix@xscale*(#1)},
  matrix yscale/.code=%
    \pgfmathsetmacro\pgf@matrix@yscale{\pgf@matrix@yscale*(#1)},
  matrix scale/.style={/tikz/matrix xscale={#1},/tikz/matrix yscale={#1}}}
\def\pgf@matrix@xscale{1}
\def\pgf@matrix@yscale{1}
\makeatother
 
\begin{document}

\title[Centralizers associated with minimal nilpotents in $F(4)$ and $G(3)$]{Centralizers associated with minimal nilpotent elements in the exceptional Lie superalgebras $F(4)$ and $G(3)$}

\author[M. Seo]{Myeongin Seo}
\address[M. Seo]{Department of Mathematical Sciences, Seoul National University, Gwanak-ro 1, Gwanak-gu, Seoul 08826, Korea}
\email{myeonginseo25@snu.ac.kr}

\date{} 

\begin{abstract}
The even part of each of the exceptional Lie superalgebras $F(4)$ and $G(3)$ has two simple summands.
For each summand, we consider the $\mathfrak{sl}_2$-triples associated with its minimal nonzero nilpotent orbit and study their extensions to $\mathfrak{osp}(1\vert 2)$-subalgebras of the full Lie superalgebra.
We give constructive proofs that these extensions exist and that those arising from the same summand form a single conjugacy class under inner automorphisms.
Thus each Lie superalgebra has exactly two such classes, one for each simple summand.
For a representative of each class, we compute the centralizers of the $\mathfrak{osp}(1\vert 2)$-subalgebra and its even $\mathfrak{sl}_2$-subalgebra explicitly.
We determine the isomorphism type of the former centralizer as a Lie superalgebra, construct a complementary submodule inside the latter, and determine its module isomorphism type.
\end{abstract}

\maketitle

\section{Introduction}

Let $\mathfrak a$ be a finite-dimensional complex semisimple Lie algebra and $A$ its adjoint group. The nilpotent elements of $\mathfrak a$ form finitely many $A$-orbits, classified through methods including weighted Dynkin diagrams and Bala--Carter theory \cites{Dynkin,BC1,BC2,CM}. The Jacobson--Morozov theorem and Kostant's conjugacy theorem connect nonzero nilpotent orbits with conjugacy classes of $\mathfrak{sl}_2$-triples \cites{Morozov,Jacobson,Kostant}. Among these orbits, the unique principal orbit is open and dense in the nilpotent cone \cite{CM}*{Theorem 4.1.6}. If $\mathfrak a$ is simple, it also has a unique nonzero nilpotent orbit of minimal dimension, which is the $A$-orbit of a highest root vector \cite{CM}*{Theorem 4.3.3}. An $\mathfrak{sl}_2$-triple containing a highest root vector induces, through its semisimple element, a grading of $\mathfrak a$ in degrees $-2,-1,0,1,2$, with one-dimensional components in degrees $\pm2$.

For a basic simple Lie superalgebra $\mathfrak g=\mathfrak g_{\bar0}\oplus\mathfrak g_{\bar1}$, the even part is reductive, so its nilpotent orbits can still be studied using Lie algebra theory. An $\mathfrak{sl}_2$-triple in $\mathfrak g_{\bar0}$ also acts on $\mathfrak g_{\bar1}$, however, and its grading and centralizer in the full Lie superalgebra reflect this additional structure. Even nilpotent orbits have accordingly been studied in connection with representations, good gradings, and centralizers of Lie superalgebras \cites{Musson,Hoyt,Han}. The literature also investigates nilpotent loci in the odd part and their orbit geometry \cites{DS,JN}.

For each nilpotent orbit in $\mathfrak{g}_{\bar{0}}$, one can construct a $W$-algebra \cite{KRW}.
Also, a nonzero nilpotent element of $\mathfrak{g}_{\bar{1}}$ contained in an $\mathfrak{osp}(1\vert 2)$-subalgebra gives rise to a supersymmetric $W$-algebra \cites{MR,MRS}.
Of particular interest here are the minimal nilpotents.
The minimality condition is given by the minimal gradings studied by Kac and Wakimoto \cite{KW}. Such a grading comes from an even $\mathfrak{sl}_2$-triple whose action on the full Lie superalgebra has weights only from $-2$ to $2$, with one-dimensional extreme weight spaces. Thus its defining condition involves the odd part as well as the even part. In $F(4)$ and $G(3)$, the minimal nonzero nilpotent orbit of each simple summand of the even part gives a minimal grading \cite{KW}.

In their study of constrained WZW models, Frappat, Ragoucy, and Sorba describe $\mathfrak{osp}(1\vert 2)$-embeddings in $F(4)$ and $G(3)$ and the corresponding decompositions of the adjoint representations \cite{FRS}.
Entova-Aizenbud and Serganova prove an analogue of the Jacobson--Morozov theorem for quasi-reductive algebraic supergroups \cite{ES}.
This theorem gives a criterion for a nonzero odd nilpotent element to lie in the image of an $\mathfrak{osp}(1\vert 2)$-embedding and shows that two such embeddings are conjugate under the even subgroup if they agree on a nonzero odd element.

The $W$-algebra associated with a minimal even nilpotent orbit is closely related to superconformal vertex operator algebras, and its structure is well understood \cites{KRW,KW}. These superconformal vertex operator algebras are also related to the supersymmetric $W$-algebras associated with $\mathfrak{osp}(1\vert 2)$-subalgebras containing the minimal even nilpotent element \cite{LSS}. Understanding the centralizers of such an $\mathfrak{osp}(1\vert 2)$-subalgebra and its even part is very useful for studying the structure of the corresponding supersymmetric $W$-algebra \cite{LSS}.

For the minimal even $\mathfrak{sl}_2$-triples in $F(4)$ and $G(3)$, the isomorphism types of their centralizers are recorded in \cite{KRW}*{Proposition 4.1}. Han's descriptions of centralizers of even nilpotent elements in these superalgebras also include their degree-zero components, which are the centralizers of the associated triples \cite{Han}. The action of an $\mathfrak{osp}(1\vert 2)$-subalgebra's centralizer on the centralizer of its even triple has been studied in other basic Lie superalgebras \cite{LSS}, which is one of the main motivations for this paper. We study the extensions of minimal even triples in $F(4)$ and $G(3)$ up to inner automorphism, together with the resulting centralizers and their actions on the triple centralizers.

Let $\mathfrak g$ be $F(4)$ or $G(3)$, and let $\mathcal T=(X,H,Y)$ be a minimal $\mathfrak{sl}_2$-triple (see \Cref{def:minimal_sl2_triple}). When odd elements $x,y$ extend $\mathcal T$ to an $\mathfrak{osp}(1\vert 2)$-subalgebra, we call $\mathcal Q=(X,x,H,y,Y)$ a minimal $\mathfrak{osp}(1\vert 2)$-quintuple (see \Cref{def:osp_1_2_quintuple,def:minimal_osp_quintuple}). We write $\mathfrak g^{\mathcal T}$ and $\mathfrak g^{\mathcal Q}$ for their centralizers in $\mathfrak g$. We prove that every minimal triple admits such an extension and that the extensions associated with each simple summand of $\mathfrak g_{\bar0}$ form a single conjugacy class under inner automorphisms. Thus each Lie superalgebra has exactly two classes. In every case, we construct a $\mathfrak g^{\mathcal Q}$-submodule $M_{\mathcal Q}^{\mathfrak g}$ such that
\[
\mathfrak g^{\mathcal T}
=\mathfrak g^{\mathcal Q}\oplus M_{\mathcal Q}^{\mathfrak g}.
\]

\begin{rem}
  The same conclusion about the number of conjugacy classes of minimal $\mathfrak{osp}(1\vert 2)$-quintuples can be deduced from \cite{FRS}. What we provide in this paper is a constructive proof for each case.
\end{rem}

The results can be summarized by the simple summand $\mathfrak a\subseteq\mathfrak g_{\bar0}$ containing $\mathcal T$. For $\mathfrak g=F(4)$, the isomorphism types are
\[
\begin{array}{c|ccc}
\mathfrak a
&\mathfrak g^{\mathcal T}
&\mathfrak g^{\mathcal Q}
&M_{\mathcal Q}^{\mathfrak g}\\ \hline
\mathfrak{sl}_2
&\mathfrak{so}_7
&G_2
&\mathbb C^7\\
\mathfrak{so}_7
&D(2,1;2)
&\mathfrak{sl}_2\oplus\mathfrak{osp}(1\vert 2)
&\mathbb C^3\otimes\mathbb C^{1\vert 2}
\end{array}
\]
and for $\mathfrak g=G(3)$, they are
\[
\begin{array}{c|ccc}
\mathfrak a
&\mathfrak g^{\mathcal T}
&\mathfrak g^{\mathcal Q}
&M_{\mathcal Q}^{\mathfrak g}\\ \hline
\mathfrak{sl}_2
&G_2
&\mathfrak{sl}_3
&\mathbb C^3\oplus(\mathbb C^3)^*\\
G_2
&\mathfrak{osp}(3\vert 2)
&\mathfrak{osp}(1\vert 2)
&\mathbb C^{3\vert 4}
\end{array}
\]
Here $\mathbb C^7$ is the seven-dimensional irreducible $G_2$-module. In the $F(4)$ table, $\mathbb C^3$ is the three-dimensional irreducible $\mathfrak{sl}_2$-module; in the $G(3)$ table, it is the natural $\mathfrak{sl}_3$-module. The symbols $\mathbb C^{1\vert 2}$ and $\mathbb C^{3\vert 4}$ denote irreducible $\mathfrak{osp}(1\vert 2)$-modules of the indicated even and odd dimensions. The full statements appear in \Cref{thm:main_F(4),thm:main_G(3)}.

\Cref{sec:preliminaries} recalls the definitions and fixes the notation used throughout the paper. \Cref{sec:shared_structure} describes the shared structure of $F(4)$ and $G(3)$ and classifies their minimal $\mathfrak{sl}_2$-triples. \Cref{sec:F4,sec:G3} give explicit bases for $F(4)$ and $G(3)$, respectively, classify the extensions of their minimal triples to $\mathfrak{osp}(1\vert 2)$-quintuples, and determine the associated centralizers and complementary submodules. \Cref{sec:commutation_relations} records all nonzero commutation relations, up to skew-supersymmetry, among the elements of the chosen bases of $F(4)$ and $G(3)$.

\section{Preliminaries}\label{sec:preliminaries}

Throughout this paper, all vector superspaces are defined over the field $\mathbb{C}$ of complex numbers.
Let $\mathbb{Z}_2=\{\bar{0},\bar{1}\}$ be the group of two elements.

\subsection{General Definitions}

We recall the notion of Lie superalgebras.
In this and the next subsection, we closely follow a standard text on Lie superalgebras \cite{CW}*{Chapter 1}.

\begin{defn}
  A \textbf{vector superspace} is a vector space $V$ endowed with a $\mathbb{Z}_2$-gradation: $V=V_{\bar{0}}\oplus V_{\bar{1}}$.
  Its \textbf{dimension} is the tuple $\dim V=(\dim V_{\bar{0}}\vert\dim V_{\bar{1}})$.
  We call $V$ \textbf{finite-dimensional} if its underlying vector space is finite-dimensional.
  We denote by $\mathbb{C}^{m\vert n}$ the vector superspace with even subspace $\mathbb{C}^m$ and odd subspace $\mathbb{C}^n$.
  The \textbf{parity} of a nonzero homogeneous element $a\in V_i$ is denoted by $|a|=i$, $i\in\mathbb{Z}_2$.
  An element in $V_{\bar{0}}$ is called \textbf{even}, while an element in $V_{\bar{1}}$ is called \textbf{odd}.
  In particular, $0$ is both even and odd.
  A \textbf{subspace} of a vector superspace $V=V_{\bar{0}}\oplus V_{\bar{1}}$ is a usual subspace $W$ of $V$ with $W=(W\cap V_{\bar{0}})\oplus(W\cap V_{\bar{1}})$.
  Here, we write $W_{\bar{0}}=W\cap V_{\bar{0}}$ and $W_{\bar{1}}=W\cap V_{\bar{1}}$.
\end{defn}

A \textbf{superalgebra} is a vector superspace $A=A_{\bar{0}}\oplus A_{\bar{1}}$ equipped with a bilinear multiplication satisfying $A_iA_j\subseteq A_{i+j}$ for all $i,j\in\mathbb{Z}_2$.
By a \textbf{subalgebra} (resp. an \textbf{ideal}), we mean a usual subalgebra (resp. ideal) that is also a subspace in the sense of vector superspaces.

\begin{defn}
  A \textbf{Lie superalgebra} is a superalgebra $\mathfrak{g}=\mathfrak{g}_{\bar{0}}\oplus\mathfrak{g}_{\bar{1}}$, whose bilinear multiplication is denoted by $[\cdot,\cdot]$, satisfying the following two additional axioms: for nonzero homogeneous elements $a,b,c\in\mathfrak{g}$,
  \begin{enumerate}
    \item (Skew-supersymmetry) $[a,b]=-(-1)^{|a||b|}[b,a]$.
    \item (Super Jacobi identity) $[a,[b,c]]=[[a,b],c]+(-1)^{|a||b|}[b,[a,c]]$.
  \end{enumerate}
\end{defn}

For a Lie superalgebra $\mathfrak{g}=\mathfrak{g}_{\bar{0}}\oplus\mathfrak{g}_{\bar{1}}$, the even part $\mathfrak{g}_{\bar{0}}$ is a Lie algebra, and the bracket defines a $\mathfrak{g}_{\bar{0}}$-module structure on $\mathfrak{g}_{\bar{1}}$ by
\[a\cdot v=[a,v]\qquad(a\in\mathfrak{g}_{\bar{0}},\ v\in\mathfrak{g}_{\bar{1}}).\]

\begin{ex}\label{ex:general_linear_lie_superalgebra}
  Let $V=V_{\bar{0}}\oplus V_{\bar{1}}$ be a vector superspace.
  The endomorphism algebra $\operatorname{End}_{\mathbb{C}}(V)$ becomes a vector superspace by setting
  \[
    \operatorname{End}_{\mathbb{C}}(V)_i
    =\{f\in\operatorname{End}_{\mathbb{C}}(V)\mid f(V_j)\subseteq V_{i+j}\text{ for all }j\in\mathbb{Z}_2\},
    \qquad i\in\mathbb{Z}_2.
  \]
  Together with the bracket $[f,g]=fg-(-1)^{|f||g|}gf$ for nonzero homogeneous endomorphisms $f$ and $g$, this is the \textbf{general linear Lie superalgebra} $\mathfrak{gl}(V)$.
  When $V=\mathbb{C}^{m\vert n}$, we write $\mathfrak{gl}(m\vert n)$ for $\mathfrak{gl}(V)$.
\end{ex}

\begin{ex}\label{ex:osp_1_2}
  The Lie superalgebra $\mathfrak{osp}(1\vert 2)$ has dimension $(3\vert 2)$, with even basis $\{X,H,Y\}$ and odd basis $\{x,y\}$.
  Its bracket is determined by skew-supersymmetry and the relations
  \begin{equation}\label{eq:osp_1_2_relations}
    \begin{gathered}
      {}[H,X]=2X,\qquad [H,Y]=-2Y,\qquad [X,Y]=H,\\
      [H,x]=x,\qquad [H,y]=-y,\qquad [X,x]=[Y,y]=0,\\
      [X,y]=x,\qquad [Y,x]=y,\\
      [x,x]=2X,\qquad [x,y]=-H,\qquad [y,y]=-2Y.
    \end{gathered}
  \end{equation}
\end{ex}

A Lie superalgebra $\mathfrak{g}$ is \textbf{abelian} if $[\mathfrak{g},\mathfrak{g}]=0$.
A nonabelian Lie superalgebra that has no nonzero proper ideal is called \textbf{simple}.
The following criterion for simplicity is used in \Cref{sec:F4,sec:G3}.

\begin{lem}\label{lem:simplicity_from_odd_part}
  Let $\mathfrak{g}=\mathfrak{g}_{\bar{0}}\oplus\mathfrak{g}_{\bar{1}}$ be a Lie superalgebra such that
  \[[\mathfrak{g}_{\bar{1}},\mathfrak{g}_{\bar{1}}]=\mathfrak{g}_{\bar{0}}\neq0.\]
  If $\mathfrak{g}_{\bar{1}}$ is a faithful irreducible $\mathfrak{g}_{\bar{0}}$-module, then $\mathfrak{g}$ is simple.
\end{lem}

\begin{proof}
  The bracket hypothesis implies that $\mathfrak{g}$ is nonabelian.
  Let $\mathfrak{i}=\mathfrak{i}_{\bar{0}}\oplus\mathfrak{i}_{\bar{1}}$ be an ideal of $\mathfrak{g}$.
  If $\mathfrak{i}_{\bar{1}}=0$, then $[\mathfrak{i}_{\bar{0}},\mathfrak{g}_{\bar{1}}]=0$, so faithfulness gives $\mathfrak{i}_{\bar{0}}=0$ and hence $\mathfrak{i}=0$.
  If $\mathfrak{i}_{\bar{1}}\neq0$, then irreducibility gives $\mathfrak{i}_{\bar{1}}=\mathfrak{g}_{\bar{1}}$.
  It follows that
  \[\mathfrak{g}_{\bar{0}}=[\mathfrak{g}_{\bar{1}},\mathfrak{g}_{\bar{1}}]=[\mathfrak{i}_{\bar{1}},\mathfrak{g}_{\bar{1}}]\subseteq\mathfrak{i}_{\bar{0}},\]
  and hence $\mathfrak{i}=\mathfrak{g}$.
\end{proof}

\begin{defn}\label{def:lie_superalgebra_homomorphism}
  Let $\mathfrak{g}$ and $\mathfrak{h}$ be Lie superalgebras.
  A \textbf{Lie superalgebra homomorphism} from $\mathfrak{g}$ to $\mathfrak{h}$ is a linear map $\varphi:\mathfrak{g}\to\mathfrak{h}$ such that
  \[
    \varphi(\mathfrak{g}_i)\subseteq\mathfrak{h}_i\quad(i\in\mathbb{Z}_2)
    \qquad\text{and}\qquad
    \varphi([a,b])=[\varphi(a),\varphi(b)]\quad(a,b\in\mathfrak{g}).
  \]
  A bijective Lie superalgebra homomorphism is called an \textbf{isomorphism}.
  An isomorphism from $\mathfrak{g}$ to itself is called an \textbf{automorphism} of $\mathfrak{g}$.
\end{defn}

\begin{ex}\label{ex:adjoint_map}
  Let $\mathfrak{g}$ be a Lie superalgebra.
  The \textbf{adjoint map} $\mathrm{ad}:\mathfrak{g}\to\mathfrak{gl}(\mathfrak{g})$ is defined by \[\mathrm{ad}(a)(b)=[a,b]\quad(a,b\in\mathfrak{g}).\]
  It follows from the super Jacobi identity that the adjoint map is a Lie superalgebra homomorphism.
\end{ex}

More generally, a \textbf{module} over a Lie superalgebra $\mathfrak{g}$ is a vector superspace $M=M_{\bar{0}}\oplus M_{\bar{1}}$ with a bilinear action $\mathfrak{g}\times M\to M$, written $(a,m)\mapsto a\cdot m$, such that $\mathfrak{g}_i\cdot M_j\subseteq M_{i+j}$ for all $i,j\in\mathbb{Z}_2$ and
\[
  [a,b]\cdot m=a\cdot(b\cdot m)-(-1)^{|a||b|}b\cdot(a\cdot m)
\]
for nonzero homogeneous $a,b\in\mathfrak{g}$ and $m\in M$.
Such an action gives a Lie superalgebra homomorphism $\rho_M:\mathfrak{g}\to\mathfrak{gl}(M)$ by $\rho_M(a)(m)=a\cdot m$.

\begin{ex}
  Let $\mathfrak{g}$ be a Lie superalgebra.
  If $a\in\mathfrak{g}$ is ad-nilpotent (that is, $\mathrm{ad}(a)$ is nilpotent), then we can consider the well-defined exponential \[\exp(\mathrm{ad}(a))=\sum_{n=0}^{\infty}\frac{1}{n!}(\mathrm{ad}(a))^n.\]
  This exponential is a Lie superalgebra automorphism if $a$ is even.
  A finite composition of $\exp(\mathrm{ad}(a))$ for even ad-nilpotent elements $a$ is called an \textbf{inner automorphism}.
\end{ex}

\begin{defn}
  Two elements $a,b\in\mathfrak{g}$ are \textbf{conjugate} if there is an inner automorphism $\sigma$ of $\mathfrak{g}$ such that $b=\sigma(a)$.
  Similarly, ordered tuples $(a_1,\dotsc,a_r)$ and $(b_1,\dotsc,b_r)$ of elements of $\mathfrak{g}$ are \textbf{conjugate} if there is an inner automorphism $\sigma$ of $\mathfrak{g}$ such that $b_i=\sigma(a_i)$ for all $i$.
\end{defn}

\subsection{Basic Simple Lie Superalgebras and Roots}

The finite-dimensional simple Lie superalgebras over $\mathbb{C}$ were classified by Kac \cite{Kac}.
Most of the Lie superalgebras considered in this paper belong to the class of basic simple Lie superalgebras, whose members are listed below; see also \cite{CW}*{Definition 1.14}.

\begin{defn}\label{def:basic_simple_lie_superalgebra}
  A finite-dimensional simple Lie superalgebra is called \textbf{basic} if it is isomorphic to one of the following Lie superalgebras:
  \begin{enumerate}
    \item[(1)] The classical series: $\mathfrak{sl}(m\vert n)$, for $m>n\geq 1$ (excluding $(m,n)=(2,1)$); $\mathfrak{sl}(m\vert m)/\mathbb{C}I_{m\vert m}$, for $m\geq 3$; and $\mathfrak{spo}(2m\vert n)$, for $m,n\geq 1$.
    \item[(2)] The exceptional cases: $D(2,1;\alpha)$, for $\alpha\in(\mathbb{C}\setminus\{0,\pm 1,-2,-\frac{1}{2}\})/\mathfrak{S}_3$; $F(4)$; and $G(3)$.
  \end{enumerate}
\end{defn}

The even part $\mathfrak{g}_{\bar{0}}$ of a basic simple Lie superalgebra $\mathfrak{g}$ is a reductive Lie algebra.
We define a \textbf{Cartan subalgebra} of $\mathfrak{g}$ as a Cartan subalgebra $\mathfrak{h}$ of $\mathfrak{g}_{\bar{0}}$.

\begin{defn}\label{def:fundamental_system}
  Let $\mathfrak{g}$ be a basic simple Lie superalgebra and $\mathfrak{h}$ a Cartan subalgebra of $\mathfrak{g}$.
  For $\alpha\in\mathfrak{h}^*$, let
  \[\mathfrak{g}_{\alpha}=\{g\in\mathfrak{g}\mid[h,g]=\alpha(h)g\text{ for all }h\in\mathfrak{h}\}.\]
  The \textbf{root system} of $\mathfrak{g}$ relative to $\mathfrak{h}$ is
  \[\Phi=\{\alpha\in\mathfrak{h}^*\setminus\{0\}\mid\mathfrak{g}_{\alpha}\neq 0\}.\]
  A subset $\Pi$ of $\Phi$ is called a \textbf{fundamental system} if it is linearly independent and every $\beta\in\Phi$ can be written as
  \[\beta=\sum_{\alpha\in\Pi}n_{\alpha}\alpha,\qquad n_{\alpha}\in\mathbb{Z},\]
  with either all $n_{\alpha}\geq0$ or all $n_{\alpha}\leq0$.
  The \textbf{positive system} determined by $\Pi$ is
  \[\Phi^+=\left\{\beta\in\Phi\,\middle|\,\beta=\sum_{\alpha\in\Pi}n_{\alpha}\alpha\text{ with all }n_{\alpha}\geq0\right\}.\]
  The elements of $\Phi$, $\Pi$, and $\Phi^+$ are called \textbf{roots}, \textbf{simple roots}, and \textbf{positive roots}, respectively.
\end{defn}

A root $\alpha$ is called \textbf{even} if $\mathfrak{g}_{\alpha}\cap\mathfrak{g}_{\bar{0}}\neq0$ and \textbf{odd} if $\mathfrak{g}_{\alpha}\cap\mathfrak{g}_{\bar{1}}\neq0$.
These two conditions need not be mutually exclusive for an arbitrary Lie superalgebra.
For a basic simple Lie superalgebra, however, the standard structure theorem gives
\[\mathfrak{g}=\mathfrak{h}\oplus\bigoplus_{\alpha\in\Phi}\mathfrak{g}_{\alpha},\]
and every nonzero root space $\mathfrak{g}_{\alpha}$ has dimension either $(1\vert0)$ or $(0\vert1)$ \cite{CW}*{Theorem 1.18}.
Consequently, every root of a basic simple Lie superalgebra is either even or odd, but not both.

\subsection{Minimal \texorpdfstring{$\mathfrak{sl}_2$}{sl2}-Triples and \texorpdfstring{$\mathfrak{osp}(1\vert 2)$}{osp(1|2)}-Quintuples}\label{subsec:minimal_nilpotents}

\begin{defn}
  Let $\mathfrak{g}$ be a Lie superalgebra.
  An $\bm{\mathfrak{sl}_2}$\textbf{-triple} in $\mathfrak{g}$ is a triple $(X,H,Y)$ of nonzero even elements of $\mathfrak{g}$ satisfying the following relations:
  \begin{equation}\label{eq:sl_2}
    [H,X]=2X,\qquad[H,Y]=-2Y,\qquad[X,Y]=H.
  \end{equation}
\end{defn}

If $(X,H,Y)$ is an $\mathfrak{sl}_2$-triple in $\mathfrak{g}$, then $\operatorname{span}_{\mathbb{C}}\{X,H,Y\}$ is an $\mathfrak{sl}_2$-subalgebra of $\mathfrak{g}$.
Thus, the adjoint action makes $\mathfrak{g}$ an $\mathfrak{sl}_2$-module, and $\mathfrak{g}$ admits the weight space decomposition
\begin{equation}\label{eq:sl2_weight_space_decomposition}
  \mathfrak{g}=\bigoplus_{n\in\mathbb{Z}}\mathfrak{g}_n.
\end{equation}

\begin{defn}\label{def:minimal_sl2_triple}
  Let $\mathfrak{g}$ be a basic simple Lie superalgebra.
  An $\mathfrak{sl}_2$-triple $(X,H,Y)$ in $\mathfrak{g}$ is \textbf{minimal} if its associated weight space decomposition \eqref{eq:sl2_weight_space_decomposition} satisfies the following:
  \begin{enumerate}
    \item[(1)] $\mathfrak{g}_n=0$ whenever $|n|>2$.
    \item[(2)] $\mathfrak{g}_2=\operatorname{span}_{\mathbb{C}}\{X\}$ and $\mathfrak{g}_{-2}=\operatorname{span}_{\mathbb{C}}\{Y\}$.
  \end{enumerate}
\end{defn}

\begin{rem}\label{rem:scalar_multiple_minimal_triple}
  If $(X,H,Y)$ is a minimal $\mathfrak{sl}_2$-triple and $c\in\mathbb{C}^{\times}$, then $(X/c,H,cY)$ is also minimal.
  Indeed, it is an $\mathfrak{sl}_2$-triple with the same weight space decomposition and extremal weight spaces.
\end{rem}

\begin{defn}\label{def:osp_1_2_quintuple}
  Let $\mathfrak{g}$ be a Lie superalgebra.
  An $\bm{\mathfrak{osp}(1\vert 2)}$\textbf{-quintuple} in $\mathfrak{g}$ is a quintuple $(X,x,H,y,Y)$ of nonzero elements such that $X,H,Y\in\mathfrak{g}_{\bar{0}}$, $x,y\in\mathfrak{g}_{\bar{1}}$, and the relations in \eqref{eq:osp_1_2_relations} hold.
  The subspace $\operatorname{span}_{\mathbb{C}}\{X,x,H,y,Y\}$ is then a Lie subalgebra of $\mathfrak{g}$ isomorphic to $\mathfrak{osp}(1\vert 2)$.
\end{defn}

\begin{defn}\label{def:minimal_osp_quintuple}
  Let $\mathfrak{g}$ be a basic simple Lie superalgebra.
  An $\mathfrak{osp}(1\vert 2)$-quintuple $(X,x,H,y,Y)$ in $\mathfrak{g}$ is \textbf{minimal} if its $\mathfrak{sl}_2$-triple $(X,H,Y)$ is minimal.
\end{defn}

The following criterion allows us to complete a minimal $\mathfrak{sl}_2$-triple to a minimal $\mathfrak{osp}(1\vert 2)$-quintuple without checking all the relations directly.

\begin{lem}\label{lem:minimality_criterion}
  Let $(X,H,Y)$ be a minimal $\mathfrak{sl}_2$-triple in a basic simple Lie superalgebra $\mathfrak{g}$.
  If $y$ is an odd element of $\mathfrak{g}$ with \[[y,y]=-2Y\qquad\text{and}\qquad[H,y]=-y,\] then $(X,[X,y],H,y,Y)$ is a minimal $\mathfrak{osp}(1\vert 2)$-quintuple.
\end{lem}

\begin{proof}
  Since $y$ is odd and $[y,y]$ is even, the skew-supersymmetry and super Jacobi identity give
  \[
    [y,[y,y]]
      =\frac{1}{3}\bigl([y,[y,y]]-[[y,y],y]+[y,[y,y]]\bigr)
      =0.
  \]
  Since $[y,y]=-2Y$, it follows that
  \[
    [Y,y]=\frac{1}{2}[y,[y,y]]=0.
  \]
  Together with $[H,y]=-y$, this shows that $y$ is a lowest weight vector of weight $-1$ for the $\mathfrak{sl}_2$-action.
  Let $x=[X,y]$.
  Finite-dimensional $\mathfrak{sl}_2$-representation theory shows that the submodule generated by $y$ is two-dimensional; in particular, $[H,x]=x$ and $[X,x]=0$.
  Moreover,
  \[
    [x,y]
      =\frac{1}{2}\bigl([x,y]+[y,x]\bigr)
      =\frac{1}{2}\bigl([[X,y],y]+[y,[X,y]]\bigr)
      =\frac{1}{2}[X,[y,y]]
      =-[X,Y]
      =-H,
  \]
  so $x$ is nonzero, and
  \[
    [Y,x]=[Y,[X,y]]=[-H,y]=y,
    \qquad
    [x,x]=[x,[X,y]]=[X,[x,y]]=2X.
  \]
  Thus, $(X,x,H,y,Y)$ is a minimal $\mathfrak{osp}(1\vert 2)$-quintuple.
  This completes the proof.
\end{proof}

We now introduce notation for the centralizers that form the main objects of this paper and for their complementary submodules.

\begin{notation}\label{not:centralizers_and_complementary_submodules}
  Let $\mathcal{Q}=(X,x,H,y,Y)$ be a minimal $\mathfrak{osp}(1\vert 2)$-quintuple in a basic simple Lie superalgebra $\mathfrak{g}$, and let $\mathcal{T}=(X,H,Y)$ be its $\mathfrak{sl}_2$-triple.
  Write $\mathfrak{g}^{\mathcal{T}}$ and $\mathfrak{g}^{\mathcal{Q}}$ for the joint centralizers in $\mathfrak{g}$ of the entries of $\mathcal{T}$ and $\mathcal{Q}$, respectively.
  Since $H=[X,Y]$ and $x=[X,y]$, these centralizers can be computed as
  \[
    \mathfrak{g}^{\mathcal{T}}=\{a\in\mathfrak{g}\mid[a,X]=[a,Y]=0\},
    \qquad
    \mathfrak{g}^{\mathcal{Q}}=\{a\in\mathfrak{g}^{\mathcal{T}}\mid[a,y]=0\}.
  \]
  If $\mathfrak{g}=\bigoplus_{n\in\mathbb{Z}}\mathfrak{g}_n$ is the weight space decomposition associated with $\mathcal{T}$, then these centralizers can also be written as
  \[
    \mathfrak{g}^{\mathcal{T}}=\{a\in\mathfrak{g}_0\mid[a,Y]=0\}=\mathfrak{g}_0^Y,
    \qquad
    \mathfrak{g}^{\mathcal{Q}}=\{a\in\mathfrak{g}_0\mid[a,y]=0\}=\mathfrak{g}_0^y.
  \]
  A weight-zero element commuting with $Y$ also commutes with $X$ by $\mathfrak{sl}_2$-representation theory, and an element commuting with $y$ also commutes with $Y=-[y,y]/2$.
  The expressions $\mathfrak{g}_0^Y$ and $\mathfrak{g}_0^y$ correspond to the notation in \cite{LSS}, where $Y$ and $y$ are written $F$ and $f$.
  The Lie superalgebra $\mathfrak{g}^{\mathcal{T}}$ contains $\mathfrak{g}^{\mathcal{Q}}$ as a Lie subalgebra, and the adjoint action makes $\mathfrak{g}^{\mathcal{T}}$ a $\mathfrak{g}^{\mathcal{Q}}$-module.
  If there exists a $\mathfrak{g}^{\mathcal{Q}}$-submodule of $\mathfrak{g}^{\mathcal{T}}$ complementary to $\mathfrak{g}^{\mathcal{Q}}$, we denote one such submodule by $M_{\mathcal{Q}}^{\mathfrak{g}}$; thus
  \[\mathfrak{g}^{\mathcal{T}}=\mathfrak{g}^{\mathcal{Q}}\oplus M_{\mathcal{Q}}^{\mathfrak{g}}.\]
\end{notation}

These objects are compatible with conjugacy.
If two minimal $\mathfrak{osp}(1\vert 2)$-quintuples $\mathcal{Q}$ and $\mathcal{Q}'$, with even triples $\mathcal{T}$ and $\mathcal{T}'$, are conjugate via an inner automorphism $\sigma$, then $\sigma$ sends $\mathfrak{g}^{\mathcal{T}}$ and $\mathfrak{g}^{\mathcal{Q}}$ isomorphically onto $\mathfrak{g}^{\mathcal{T}'}$ and $\mathfrak{g}^{\mathcal{Q}'}$, respectively, and transports any chosen complementary submodule accordingly.
Thus, it suffices to compute these objects for one representative of each conjugacy class.

\section{Shared Structure and Minimal \texorpdfstring{$\mathfrak{sl}_2$}{sl2}-Triples in \texorpdfstring{$F(4)$}{F(4)} and \texorpdfstring{$G(3)$}{G(3)}}\label{sec:shared_structure}

\subsection{Shared Structure of \texorpdfstring{$F(4)$}{F(4)} and \texorpdfstring{$G(3)$}{G(3)}}\label{subsec:structure_of_F4_and_G3}

In this subsection, we fix the structural data and basis conventions shared by the basic simple Lie superalgebras $F(4)$ and $G(3)$.
Let $\mathfrak{g}$ be either $F(4)$ or $G(3)$.
Its even part is the direct sum of two simple ideals $\mathfrak{a}_1$ and $\mathfrak{a}_2$, and its odd part, regarded as a $\mathfrak{g}_{\bar{0}}$-module, has the form
\[\mathfrak{g}_{\bar{0}}=\mathfrak{a}_1\oplus\mathfrak{a}_2,
  \qquad \mathfrak{g}_{\bar{1}}\cong V_1\otimes V_2,\]
where $\mathfrak{a}_1\cong\mathfrak{sl}_2$, $V_1$ is its $2$-dimensional irreducible module, and $V_2$ is an irreducible $\mathfrak{a}_2$-module.
The second ideal is
\[
  \mathfrak{a}_2\cong
  \begin{cases}
    \mathfrak{so}_7,&\mathfrak{g}=F(4),\\
    G_2,&\mathfrak{g}=G(3).
  \end{cases}
\]

\begin{lem}\label{lem:automorphisms_preserve_even_ideals}
  Every automorphism of $\mathfrak{g}$ preserves both $\mathfrak{a}_1$ and $\mathfrak{a}_2$.
\end{lem}

\begin{proof}
  By \Cref{def:lie_superalgebra_homomorphism}, an automorphism $\varphi$ of $\mathfrak{g}$ preserves the $\mathbb{Z}_2$-grading and hence restricts to an automorphism of $\mathfrak{g}_{\bar{0}}$.
  This restriction permutes the simple ideals of $\mathfrak{g}_{\bar{0}}$.
  Since $\mathfrak{a}_1$ and $\mathfrak{a}_2$ are not isomorphic, it follows that $\varphi(\mathfrak{a}_i)=\mathfrak{a}_i$ for $i=1,2$.
\end{proof}

Fix Cartan subalgebras $\mathfrak{h}_1$ and $\mathfrak{h}_2$ of $\mathfrak{a}_1$ and $\mathfrak{a}_2$, respectively, and set $\mathfrak{h}:=\mathfrak{h}_1\oplus\mathfrak{h}_2\subseteq\mathfrak{g}_{\bar{0}}$.
Then $\mathfrak{h}$ is a Cartan subalgebra of $\mathfrak{g}$.
We regard $\mathfrak{h}_1^*$ and $\mathfrak{h}_2^*$ as subspaces of $\mathfrak{h}^*$ by extending functionals by zero on the other summand.
Let $\Phi_i\subseteq\mathfrak{h}_i^*$ be the root system of $\mathfrak{a}_i$.
We choose a fundamental system $\Pi_1=\{\alpha\}$ of $\Phi_1$ and a fundamental system
\[\Pi_2=\{\alpha_1,\ldots,\alpha_r\}\]
of $\Phi_2$.
For $F(4)$, we take $r=3$, with $\alpha_1$ and $\alpha_2$ long and $\alpha_3$ short; for $G(3)$, we take $r=2$, with $\alpha_1$ short and $\alpha_2$ long.
Write $\Phi_i^+$ for the positive systems determined by $\Pi_i$.
For explicit constructions of these root systems, see \cite{Humphreys}*{Section 12}.

Choose nonzero root vectors $e_{\alpha}\in(\mathfrak{a}_1)_{\alpha}$ and $e_{-\alpha}\in(\mathfrak{a}_1)_{-\alpha}$, normalized so that, with $h:=[e_{\alpha},e_{-\alpha}]$, the triple $(e_{\alpha},h,e_{-\alpha})$ is an $\mathfrak{sl}_2$-triple.
For $\beta\in\Phi_i^+$ ($i=1,2$), let $h_{\beta}$ be the coroot of $\beta$, characterized by
\[h_{\beta}\in[(\mathfrak{a}_i)_{\beta},(\mathfrak{a}_i)_{-\beta}]
  \qquad\text{and}\qquad\beta(h_{\beta})=2.\]
In particular, $h_{\alpha}=h$.
For each $\beta\in\Phi_2^+$, once a nonzero root vector $e_{\beta}\in(\mathfrak{a}_2)_{\beta}$ has been chosen, let $e_{-\beta}\in(\mathfrak{a}_2)_{-\beta}$ be the unique root vector satisfying
\[[e_{\beta},e_{-\beta}]=h_{\beta}.\]
Write $h_i:=h_{\alpha_i}$ for $1\leq i\leq r$.

We use the following notation for highest weight modules throughout the remainder of the paper.

\begin{notation}\label{not:highest_weight_modules}
  Let $\mathfrak{a}$ be a finite-dimensional semisimple Lie algebra with a fixed Cartan subalgebra and fundamental system.
  For a weight $\mu$, denote by $V_{\mathfrak{a}}(\mu)$ the irreducible highest weight $\mathfrak{a}$-module with highest weight $\mu$.
\end{notation}

Let $\lambda\in\mathfrak{h}_1^*$ and $\lambda_1,\ldots,\lambda_r\in\mathfrak{h}_2^*$ be the fundamental weights relative to $\Pi_1$ and $\Pi_2$, respectively; thus
\[\lambda(h)=1,\qquad \lambda_i(h_j)=\delta_{ij}.\]
With \Cref{not:highest_weight_modules},
\[
  V_1\cong V_{\mathfrak{a}_1}(\lambda),\qquad
  V_2\cong
  \begin{cases}
    V_{\mathfrak{a}_2}(\lambda_3),&\mathfrak{g}=F(4),\\
    V_{\mathfrak{a}_2}(\lambda_1),&\mathfrak{g}=G(3).
  \end{cases}
\]
The two weights of $V_1$ are $\lambda$ and $-\lambda$.
All the weight spaces of $V_2$ are one-dimensional in both cases.
We identify $\mathfrak{g}_{\bar{1}}$ with $V_1\otimes V_2$.

For any $\mathfrak{h}$-module $M$ and $\mu\in\mathfrak{h}^*$, denote by $M_{\mu}$ the weight space $\{m\in M\mid t\cdot m=\mu(t)m\text{ for all }t\in\mathfrak{h}\}$.
The same notation applies to modules over $\mathfrak{h}_1$ and $\mathfrak{h}_2$.
Choose a nonzero vector $v_-\in(V_1)_{-\lambda}$ and set $v_+:=e_{\alpha}\cdot v_-$.
Whenever nonzero vectors $v_{\mu}\in(V_2)_{\mu}$ have been chosen, set
\[v_{\mu}^{+}:=v_+\otimes v_{\mu},\qquad v_{\mu}^{-}:=v_-\otimes v_{\mu}.\]

We shall repeatedly use two standard root string facts.
First, if $\beta$, $\gamma$, and $\beta+\gamma$ are roots of $\mathfrak{a}_2$, then
\[[(\mathfrak{a}_2)_{\beta},(\mathfrak{a}_2)_{\gamma}]=(\mathfrak{a}_2)_{\beta+\gamma}.\]
Second, if $\mu$ and $\mu+\beta$ are consecutive weights in a $\beta$-string of $V_2$, then every nonzero $e_{\beta}\in(\mathfrak{a}_2)_{\beta}$ acts nontrivially from $(V_2)_{\mu}$ to $(V_2)_{\mu+\beta}$.
Thus all vectors defined recursively in \Cref{subsec:explicit_basis_for_F(4),subsec:explicit_basis_for_G(3)} are nonzero.

\begin{prop}\label{prop:common_basis_construction}
  Suppose that nonzero vectors $e_{\beta}\in(\mathfrak{a}_2)_{\beta}$ and $v_{\mu}\in(V_2)_{\mu}$ have been chosen for every $\beta\in\Phi_2^+$ and every weight $\mu$ of $V_2$, respectively.
  Then
  \[\mathcal{B}_{\bar{0}}
    :=\{h,e_{\alpha},e_{-\alpha}\}
      \cup\{h_1,\ldots,h_r\}
      \cup\bigcup_{\beta\in\Phi_2^+}\{e_{\beta},e_{-\beta}\}
  \]
  is a basis of $\mathfrak{g}_{\bar{0}}$, and
  \[\mathcal{B}_{\bar{1}}
    :=\{v_{\mu}^{+},v_{\mu}^{-}\mid\mu\text{ is a weight of }V_2\}
  \]
  is a basis of $\mathfrak{g}_{\bar{1}}$.
  Consequently, $\mathcal{B}_{\bar{0}}\cup\mathcal{B}_{\bar{1}}$ is a basis of $\mathfrak{g}$.
\end{prop}

\begin{proof}
  The simple coroots $h_1,\ldots,h_r$ form a basis of $\mathfrak{h}_2$, while the vectors $e_{\beta}$ and $e_{-\beta}$ span the one-dimensional root spaces $(\mathfrak{a}_2)_{\beta}$ and $(\mathfrak{a}_2)_{-\beta}$, respectively.
  Hence the root space decomposition of $\mathfrak{a}_2$ shows that the elements of $\mathcal{B}_{\bar{0}}$ belonging to $\mathfrak{a}_2$ form a basis of $\mathfrak{a}_2$.
  Since $\{h,e_{\alpha},e_{-\alpha}\}$ is a basis of $\mathfrak{a}_1$, the first assertion follows.
  The vectors $v_+$ and $v_-$ form a basis of $V_1$, and the multiplicity-free property of $V_2$ implies that the chosen vectors $v_{\mu}$ form a basis of $V_2$.
  The second assertion now follows from $\mathfrak{g}_{\bar{1}}\cong V_1\otimes V_2$.
\end{proof}

\subsection{Classification of Minimal \texorpdfstring{$\mathfrak{sl}_2$}{sl2}-Triples}\label{subsec:minimal_nilpotent_orbits}

For $F(4)$ and $G(3)$, the proposition below follows from the discussion in \cite{KW}*{Section 5}.
There, $e=X/2$, $x=H/2$, and $f=Y$ in our notation, so the corresponding weight space conditions agree with those of \Cref{def:minimal_sl2_triple}.

\begin{prop}\label{prop:minimal_triple_entries}
  Let $\theta$ be the highest root in $\Phi_2^+$.
  \begin{enumerate}
    \item[(1)] The triples $(e_{\alpha},h,e_{-\alpha})$ and $(e_{\theta},h_{\theta},e_{-\theta})$ are minimal.
    \item[(2)] If $(X,H,Y)$ is a minimal $\mathfrak{sl}_2$-triple in $\mathfrak{g}$, then $Y$ belongs to either $\mathfrak{a}_1$ or $\mathfrak{a}_2$.
    In the first case, $Y$ is conjugate to $e_{-\alpha}$; in the second, to $e_{-\theta}$.
  \end{enumerate}
\end{prop}

The following conjugacy theorem of Kostant \cite{Kostant} allows us to fix the entire $\mathfrak{sl}_2$-triple, up to conjugacy, once its nilpotent entry $Y$ has been fixed.

\begin{prop}\label{prop:Kostant}
  Any two $\mathfrak{sl}_2$-triples $(X,H,Y)$ and $(X',H',Y)$ in a finite-dimensional semisimple Lie algebra that contain the same nilpotent element $Y$ are conjugate by an inner automorphism fixing $Y$.
\end{prop}

Since $\mathfrak{g}_{\bar{0}}=\mathfrak{a}_1\oplus\mathfrak{a}_2$ is semisimple, its inner automorphisms are generated by exponentials of even root vectors and therefore extend to inner automorphisms of $\mathfrak{g}$.
By \Cref{prop:minimal_triple_entries,prop:Kostant}, every minimal $\mathfrak{sl}_2$-triple is conjugate to $(e_{\alpha},h,e_{-\alpha})$ or $(e_{\theta},h_{\theta},e_{-\theta})$.
These two triples are not conjugate by \Cref{lem:automorphisms_preserve_even_ideals}, since their entries $Y$ lie in different simple ideals of $\mathfrak{g}_{\bar{0}}$.
Hence minimal $\mathfrak{sl}_2$-triples form exactly two conjugacy classes.
To classify minimal $\mathfrak{osp}(1\vert 2)$-quintuples, it therefore suffices to consider extensions of these two triples.
The defining relations and \Cref{lem:minimality_criterion} show that an odd element $y$ completes the first or second triple, respectively, precisely when
\begin{itemize}
  \item $[y,y]=-2e_{-\alpha}$ and $[h,y]=-y$; or
  \item $[y,y]=-2e_{-\theta}$ and $[h_{\theta},y]=-y$.
\end{itemize}
The case-by-case arguments in \Cref{sec:F4,sec:G3} establish the existence and uniqueness up to conjugacy of these extensions.

For the isomorphism types of the centralizers of these $\mathfrak{sl}_2$-triples, see \cite{KRW}*{Proposition 4.1}.
The case-by-case computations in \Cref{sec:F4,sec:G3} require explicit realizations of these centralizers in the chosen bases.

\subsection{Root Exponentials and Diagonal Automorphisms}\label{subsec:inner_automorphisms}

The case-by-case arguments in \Cref{sec:F4,sec:G3} use inner automorphisms both to compare odd elements while fixing a chosen $\mathfrak{sl}_2$-triple and to rescale root vectors.
We first define the root exponential automorphisms $\sigma_{\beta,c}$ and then construct the diagonal automorphisms $\tau_{\beta,t}$ from them.

Let $\beta\in\mathfrak{h}^*$ be an even root of $\mathfrak{g}$.
For $c\in\mathbb{C}$, write
\[\sigma_{\beta,c}:=\exp\bigl(\operatorname{ad}(ce_{\beta})\bigr).\]
Since $e_{\beta}$ is even and ad-nilpotent, $\sigma_{\beta,c}$ is an inner automorphism of $\mathfrak{g}$.
Also, for positive $\beta$ and $t\in\mathbb{C}^{\times}$, define
\[n_{\beta}(t):=\sigma_{\beta,t}\circ\sigma_{-\beta,-t^{-1}}\circ\sigma_{\beta,t}
  \qquad\text{and}\qquad
  \tau_{\beta,t}:=n_{\beta}(t)\circ n_{\beta}(-1).\]
Since $\beta$ is positive, the vectors $e_{\beta}$, $h_{\beta}$, and $e_{-\beta}$ span an $\mathfrak{sl}_2$-subalgebra.
In the standard two-dimensional realization of this $\mathfrak{sl}_2$-subalgebra, $e_{\beta}$, $h_{\beta}$, and $e_{-\beta}$ correspond to $E_{12}$, $\operatorname{diag}(1,-1)$, and $E_{21}$, respectively, where $E_{ij}$ denotes the usual matrix unit.
In this realization, the products defining $n_{\beta}(t)$ and $\tau_{\beta,t}$ take the following matrix form:
\[\begin{aligned}
  \begin{pmatrix}1&t\\0&1\end{pmatrix}
  \begin{pmatrix}1&0\\-t^{-1}&1\end{pmatrix}
  \begin{pmatrix}1&t\\0&1\end{pmatrix}
  &=\begin{pmatrix}0&t\\-t^{-1}&0\end{pmatrix},\\
  \begin{pmatrix}0&t\\-t^{-1}&0\end{pmatrix}
  \begin{pmatrix}0&-1\\1&0\end{pmatrix}
  &=\begin{pmatrix}t&0\\0&t^{-1}\end{pmatrix}.
\end{aligned}\]
Choose $z\in\mathbb{C}$ with $e^z=t$, so $\operatorname{diag}(t,t^{-1})=\exp\bigl(z\operatorname{diag}(1,-1)\bigr)$.
Applying the matrix identity in the corresponding $\mathrm{SL}_2(\mathbb{C})$-representation on $\mathfrak{g}$ then gives $\tau_{\beta,t}=\exp\bigl(z\,\operatorname{ad}(h_{\beta})\bigr)$.
Consequently, if $v\in\mathfrak{g}$ is nonzero and $[h_{\beta},v]=mv$, then $m\in\mathbb{Z}$ and
\[\tau_{\beta,t}(v)=e^{zm}v=t^m v.\]

The following lemma shows when the root exponentials used in these constructions fix an $\mathfrak{sl}_2$-triple pointwise.

\begin{lem}\label{lem:root_exponentials_fix_triple}
  Let $(X,H,Y)$ be an $\mathfrak{sl}_2$-triple with $H\in\mathfrak{h}$ whose associated weight space decomposition satisfies
  \[\mathfrak{g}_2=\operatorname{span}_{\mathbb{C}}\{X\}\qquad\text{and}\qquad
    \mathfrak{g}_{-2}=\operatorname{span}_{\mathbb{C}}\{Y\}.\]
  If $\beta(H)=0$, then $\sigma_{\beta,c}$ fixes $X$, $H$, and $Y$ for every $c\in\mathbb{C}$.
\end{lem}

\begin{proof}
  Since $[H,e_{\beta}]=0$, the operator $\operatorname{ad}(e_{\beta})$ preserves every $H$-weight space.
  In particular, it preserves the one-dimensional spaces $\mathfrak{g}_2$ and $\mathfrak{g}_{-2}$.
  Its restrictions to these spaces are both nilpotent and hence zero.
  Thus, $[e_{\beta},X]=[e_{\beta},H]=[e_{\beta},Y]=0$, which proves the assertion.
\end{proof}

Consequently, for a triple satisfying the hypotheses of \Cref{lem:root_exponentials_fix_triple}, every finite composition of automorphisms $\sigma_{\gamma,c}$ with $\gamma(H)=0$ fixes $(X,H,Y)$ pointwise.
For example, if $\beta(H)=0$, then $\tau_{\beta,t}$ fixes the triple, since its defining composition uses automorphisms associated with $\beta$ and $-\beta$.
In particular, these conclusions apply to a minimal $\mathfrak{sl}_2$-triple with $H\in\mathfrak{h}$ whenever the roots used vanish on $H$.

We end with an observation on the conjugacy of minimal $\mathfrak{osp}(1\vert 2)$-quintuples.
Suppose $(X,x,H,y,Y)$ and $(X,x',H,y',Y)$ are two such quintuples with $H\in\mathfrak{h}$.
Let $\varphi$ be a finite composition of automorphisms $\sigma_{\gamma,c}$ with $\gamma(H)=0$ for each root $\gamma$, and suppose $\varphi(y)=y'$.
Since $\varphi$ fixes $X$, we have
\[\varphi(x)=\varphi([X,y])=[X,\varphi(y)]=[X,y']=x'.\]
Thus $\varphi$ conjugates the two quintuples.

\section{Minimal \texorpdfstring{$\mathfrak{osp}(1\vert 2)$}{osp(1|2)}-quintuples and Their Centralizers in \texorpdfstring{$F(4)$}{F(4)}}\label{sec:F4}

Throughout this section, let $\mathfrak{g}=F(4)$.

\subsection{An Explicit Basis for \texorpdfstring{$F(4)$}{F(4)}}\label{subsec:explicit_basis_for_F(4)}

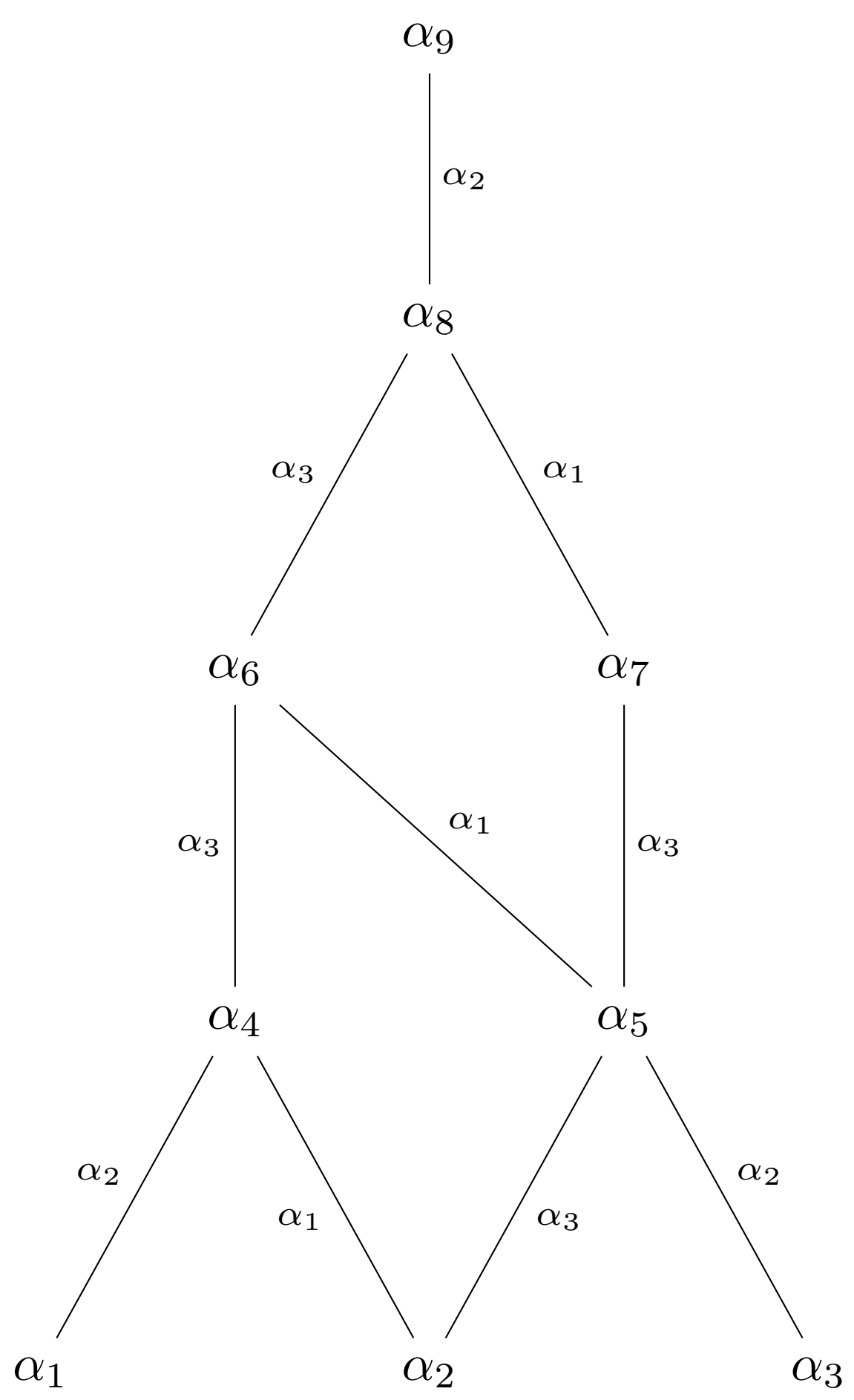
\begin{figure}[htbp]
  % https://tikzcd.yichuanshen.de/#N4Igdg9gJgpgziAXAbVABwnAlgFyxMJZAJgBoAGAXVJADcBDAGwFcYkQAdDptAC3oD6AThABfUuky58hFGQDM1Ok1bsuPfgIAcYiSAzY8BIgEZSAdiUMWbRJ26M+ggGy7JhmUXkWrK2-Y1Bczd9KSNZZG8TE18bNQcnAQBWEINpYxQzaNjVO3VHTQAWVLDPFHJSEySc-3zEkxKPDJJK6pprXICCwWJG9IjC1pr4wIF5MSUYKABzeCJQADMAJwgAWyRnGhwIJCT2vxHugV6aXhh6KHZIMDZxRZX1xHMtncQ95Ti8hM0G0-PLuzXW56ZZrJDPEDbJCDEBnC5XAhsfafLqJcY0Rj0ABGMEYAAVShkQEssNNeDgQqDHloXtC-vDAYiQMjOnVNCcQJicfjCbJiaTyZSHrtaYgyLD-gibsyPqzvoJxncQFS6ZDXuKOrV5QJfpzsbiCU0+SSyRT6QDwEylSrEDCoYhvLKtaN0XruYb+uwTYLzVLgfcwWLRWYJQzLdKWc6joqQcKHcHfYyI07DvUZVyDbyvQKKda4yH7RUU19Rhy4Ragen9TyjdnTRNREA
  \begin{tikzcd}[nodes={scale=1.3}]%[matrix scale=0.4,nodes={scale=1.3}]
                                              &                                             & \alpha_9                                                                        &                                                                                  &                                              \\
                                              &                                             &                                                                                 &                                                                                  &                                              \\
                                              &                                             &                                                                                 &                                                                                  &                                              \\
                                              &                                             & \alpha_8 \arrow[uuu, "\alpha_2"', no head]                                      &                                                                                  &                                              \\
                                              &                                             &                                                                                 &                                                                                  &                                              \\
                                              &                                             &                                                                                 &                                                                                  &                                              \\
                                              &                                             &                                                                                 &                                                                                  &                                              \\
                                              & \alpha_6 \arrow[ruuuu, "\alpha_3", no head] &                                                                                 & \alpha_7 \arrow[luuuu, "\alpha_1"', no head]                                     &                                              \\
                                              &                                             &                                                                                 &                                                                                  &                                              \\
                                              &                                             &                                                                                 &                                                                                  &                                              \\
                                              &                                             &                                                                                 &                                                                                  &                                              \\
                                              & \alpha_4 \arrow[uuuu, "\alpha_3", no head]  &                                                                                 & \alpha_5 \arrow[lluuuu, "\alpha_1"', no head] \arrow[uuuu, "\alpha_3"', no head] &                                              \\
                                              &                                             &                                                                                 &                                                                                  &                                              \\
                                              &                                             &                                                                                 &                                                                                  &                                              \\
                                              &                                             &                                                                                 &                                                                                  &                                              \\
  \alpha_1 \arrow[ruuuu, "\alpha_2", no head] &                                             & \alpha_2 \arrow[luuuu, "\alpha_1", no head] \arrow[ruuuu, "\alpha_3"', no head] &                                                                                  & \alpha_3 \arrow[luuuu, "\alpha_2"', no head]
  \end{tikzcd}
  \caption{Positive roots of $\mathfrak{so}_7$}
  \label{fig:F4-positive-roots}
\end{figure}

Retain the notation and choices of \Cref{subsec:structure_of_F4_and_G3}.
In the present case, $\mathfrak{a}_2\cong\mathfrak{so}_7$ and $\Pi_2=\{\alpha_1,\alpha_2,\alpha_3\}$, where $\alpha_1$ and $\alpha_2$ are long roots and $\alpha_3$ is a short root.
The nine positive roots of $\mathfrak{a}_2$ are shown in \Cref{fig:F4-positive-roots}.

Choose any nonzero $e_{\alpha_i}\in(\mathfrak{a}_2)_{\alpha_i}$ for $i=1,2,3$.
For $i\in\{4,5,6,7,8,9\}$, define $e_{\alpha_i}\in(\mathfrak{a}_2)_{\alpha_i}$ as follows:
\[\begin{aligned}
  &e_{\alpha_4}=[e_{\alpha_2},e_{\alpha_1}], &&\qquad e_{\alpha_5}=[e_{\alpha_2},e_{\alpha_3}], &&\qquad e_{\alpha_6}=[e_{\alpha_3},e_{\alpha_4}]=[e_{\alpha_1},e_{\alpha_5}],\\
  &e_{\alpha_7}=\frac{1}{2}[e_{\alpha_3},e_{\alpha_5}], &&\qquad e_{\alpha_8}=\frac{1}{2}[e_{\alpha_3},e_{\alpha_6}]=[e_{\alpha_1},e_{\alpha_7}], &&\qquad e_{\alpha_9}=[e_{\alpha_2},e_{\alpha_8}].
\end{aligned}\]
For $i=1,\ldots,9$, write $h_i:=h_{\alpha_i}$ and let $e_{-\alpha_i}$ be the negative root vector determined by the convention in \Cref{subsec:structure_of_F4_and_G3}.
The first assertion of \Cref{prop:common_basis_construction} gives a basis of $\mathfrak{g}_{\bar{0}}$:
\[\{h,e_{\alpha},e_{-\alpha}\}\cup\left(\{h_1,h_2,h_3\}\cup\bigcup_{i=1}^{9}\{e_{\alpha_i},e_{-\alpha_i}\}\right).\]

\begin{figure}[htbp]
  % https://tikzcd.yichuanshen.de/#N4Igdg9gJgpgziAXAbVABwnAlgFyxMJZARgBoAGAXVJADcBDAGwFcYkQAdDx+gWwCMo9APoBmEAF9S6TLnyEUZYtTpNW7LjwFDhAJgC0mvoJHipM7HgJEyulQxZtEnbsZ3FDr7SN0BqI95iktIgGJbyROSkovZqTi5aJsIeAUlmIWFy1ii60bGO7J6J7v5eacEWWQokpAAs+erORW4ixKXFPs2B6ZVW1WQArA3xXUl+qTo9obJ9NqQAbMOFE6aSKjBQAObwRKAAZgBOELxIZCA4EEhRIAAWMPRQ7JBgbDQ8-DCMAAozEc4HWE2NxwIBoDkaLiYaBuq3MIEOxyQuXOl0QZzuDyeBFeqgKzi4UJhelBIHenx+4WyIABQJBcIRJ0QohoFyRYLiGg4hNWNAxj2czxxZO+vypjBgezpIQZSFqLNRyPB8QJjGhrRJwopVXYNOBJL5WJeFXhR0ZA3lSGZuIhKrVyX19354GxGvoHxFlIUpIlUv2pqQ5pRsvZeMhqqJ4jebvJoq9upBvMdhrY9P9iHmFsQgaVnO5xKj7q1s3+gL1icxApdqcRiAA7JmM9blVzwzzSdGPdqS7SHRXnUaJJQJEA
  \begin{tikzcd}[matrix xscale=-0.4,matrix yscale=2.2,nodes={scale=1.2}]
                                                      & \lambda_3                                                                                       &                                                       \\
                                                      & \lambda_2-\lambda_3 \arrow[u, "\alpha_3"', no head]                                             &                                                       \\
                                                      & \lambda_1-\lambda_2+\lambda_3 \arrow[u, "\alpha_2"', no head]                                   &                                                       \\
  \lambda_1-\lambda_3 \arrow[ru, "\alpha_3", no head] &                                                                                                 & -\lambda_1+\lambda_3 \arrow[lu, "\alpha_1"', no head] \\
                                                      & -\lambda_1+\lambda_2-\lambda_3 \arrow[lu, "\alpha_1", no head] \arrow[ru, "\alpha_3"', no head] &                                                       \\
                                                      & -\lambda_2+\lambda_3 \arrow[u, "\alpha_2"', no head]                                            &                                                       \\
                                                      & -\lambda_3 \arrow[u, "\alpha_3"', no head]                                                      &                                                      
  \end{tikzcd}
  \caption{Weights of $V_{\mathfrak{so}_7}(\lambda_3)$}
  \label{fig:F4-weight-diagram}
\end{figure}
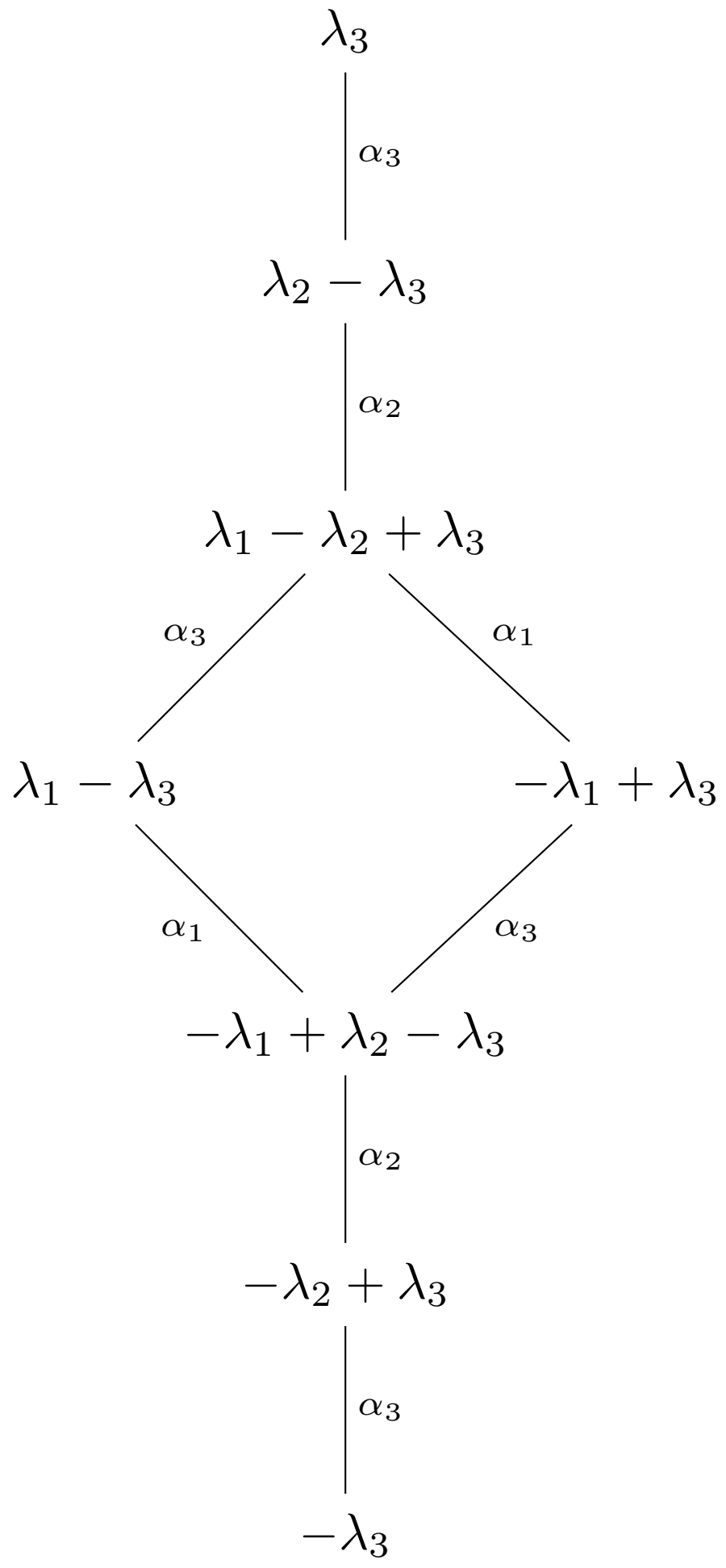

We next specify the basis vectors in $\mathfrak{g}_{\bar{1}}\cong V_1\otimes V_2$.
The weights of $V_2$ are shown in \Cref{fig:F4-weight-diagram}.

\begin{rem}
  In terms of the simple roots, the fundamental weights are given by
  \[\lambda=\frac{1}{2}\alpha,\qquad\lambda_1=\alpha_1+\alpha_2+\alpha_3,\qquad\lambda_2=\alpha_1+2\alpha_2+2\alpha_3,\qquad\lambda_3=\frac{1}{2}\alpha_1+\alpha_2+\frac{3}{2}\alpha_3.\]
\end{rem}

Choose any nonzero element $v_{-\lambda_3}\in (V_2)_{-\lambda_3}$.
In turn, for the other weights $\mu$ of $V_2$, define $v_{\mu}\in (V_2)_{\mu}$ as follows:
\[\begin{aligned}
  &v_{-\lambda_2+\lambda_3}=e_{\alpha_3}\cdot v_{-\lambda_3}, &&v_{-\lambda_1+\lambda_2-\lambda_3}=e_{\alpha_2}\cdot v_{-\lambda_2+\lambda_3}=e_{\alpha_5}\cdot v_{-\lambda_3},\\
  &v_{\lambda_1-\lambda_3}=e_{\alpha_1}\cdot v_{-\lambda_1+\lambda_2-\lambda_3}=e_{\alpha_6}\cdot v_{-\lambda_3}, &&v_{-\lambda_1+\lambda_3}=e_{\alpha_3}\cdot v_{-\lambda_1+\lambda_2-\lambda_3}=e_{\alpha_7}\cdot v_{-\lambda_3},\\
  &v_{\lambda_1-\lambda_2+\lambda_3}=e_{\alpha_3}\cdot v_{\lambda_1-\lambda_3}=e_{\alpha_1}\cdot v_{-\lambda_1+\lambda_3}=e_{\alpha_8}\cdot v_{-\lambda_3},\\
  &v_{\lambda_2-\lambda_3}=e_{\alpha_2}\cdot v_{\lambda_1-\lambda_2+\lambda_3}=e_{\alpha_9}\cdot v_{-\lambda_3}, &&v_{\lambda_3}=e_{\alpha_3}\cdot v_{\lambda_2-\lambda_3}=e_{\alpha_9}\cdot v_{-\lambda_2+\lambda_3}.
\end{aligned}\]

These vectors are nonzero by the root string facts recalled in \Cref{subsec:structure_of_F4_and_G3}.
The second assertion of \Cref{prop:common_basis_construction} therefore shows that, as $\mu$ ranges over the weights in \Cref{fig:F4-weight-diagram}, the $16$ elements $v_{\mu}^{+}$ and $v_{\mu}^{-}$ form a basis of $\mathfrak{g}_{\bar{1}}$.

\begin{rem}
  The roots $\lambda-\lambda_3$, $\alpha_3$, $\alpha_2$, and $\alpha_1$ form a fundamental system of $F(4)$:
  \begin{equation}\label{eq:F4_fundamental_system}
    \Pi=\{\lambda-\lambda_3,\alpha_3,\alpha_2,\alpha_1\}.
  \end{equation}
  Fix the order of the elements of $\Pi$ as shown in \eqref{eq:F4_fundamental_system}.
  The corresponding Cartan matrix $A=(a_{ij})$, as recorded in \cite{Dictionary}, is
  \begin{equation}\label{eq:Cartan_mat_of_F(4)}
    \begin{pmatrix}0&1&0&0\\-1&2&-2&0\\0&-1&2&-1\\0&0&-1&2\end{pmatrix}.
  \end{equation}
  Here, if we denote by $\beta_i$ the $i$-th root in $\Pi$, then $a_{ij}=\beta_j(H_i)$ for a suitably chosen nonzero $H_i\in[\mathfrak{g}_{\beta_i},\mathfrak{g}_{-\beta_i}]$.
\end{rem}

Set $\delta=\lambda-\lambda_3$.
The vectors $v_{-\lambda_3}^{+}$ and $v_{\lambda_3}^{-}$ belong to the opposite odd root spaces $\mathfrak{g}_{\delta}$ and $\mathfrak{g}_{-\delta}$, respectively.
Standard facts about root spaces in basic simple Lie superalgebras imply that their bracket is a nonzero element of $\mathfrak{h}$; see \cite{CW}*{Theorem 1.18}.
The Cartan matrix \eqref{eq:Cartan_mat_of_F(4)} shows that this bracket is a nonzero scalar multiple of
\[\frac{3}{2}h+h_1+2h_2+\frac{3}{2}h_3.\]
Replacing $v_-$ by $cv_-$ also replaces $v_+$ by $cv_+$ and multiplies the bracket by $c^2$.
Since every nonzero complex number has a square root, we may therefore normalize $v_-$ so that
\begin{equation}\label{eq:F(4)_assumption}
  \left[v_{-\lambda_3}^{+},\,v_{\lambda_3}^{-}\right]=\frac{3}{2}h+h_1+2h_2+\frac{3}{2}h_3.
\end{equation}
The complete commutation relations among these basis elements are recorded in \Cref{sec:commutation_relations_for_F(4)}.

\subsection{Classification of Minimal \texorpdfstring{$\mathfrak{osp}(1\vert 2)$}{osp(1|2)}-Quintuples in \texorpdfstring{$F(4)$}{F(4)}}\label{subsec:F4_conjugacy}

By \Cref{subsec:minimal_nilpotent_orbits}, minimal $\mathfrak{sl}_2$-triples in $F(4)$ form two distinct conjugacy classes, represented by triples associated with $\alpha$ and with the highest root $\alpha_9$ of $\mathfrak{a}_2$.
Every minimal $\mathfrak{osp}(1\vert 2)$-quintuple extends a triple in one of these classes.
Using the diagonal inner automorphisms $\tau_{\beta,t}$ constructed in \Cref{subsec:inner_automorphisms} to rescale the root vectors, we may take these triples to be
\[\left(\frac{1}{3}e_{\alpha},h,3e_{-\alpha}\right)\qquad\text{or}\qquad\left(\frac{1}{2}e_{\alpha_9},h_9,2e_{-\alpha_9}\right).\]
Note that $h_9=h_1+2h_2+h_3$.
Conjugate quintuples have conjugate even triples, so we classify the extensions of each chosen triple under automorphisms fixing that triple.

For a quintuple extending the first triple, its fourth entry $y$ satisfies $[h,y]=-y$; hence $y$ lies in the span of the $v_{\mu}^{-}$ as $\mu$ ranges over the weights of $V_2$.

\begin{prop}\label{prop:case_1_of_F(4)}
  Let
  \[y=\sum_{\mu}b_{\mu}v_{\mu}^{-}\qquad(b_{\mu}\in\mathbb{C}),\]
  where the summation runs over the weights $\mu$ of $V_2$, and set
  \begin{equation}\label{eq:F4_case_1_scalar}
    b:=b_{\lambda_1-\lambda_3}b_{-\lambda_1+\lambda_3}
      -b_{-\lambda_1+\lambda_2-\lambda_3}b_{\lambda_1-\lambda_2+\lambda_3}
      +b_{-\lambda_2+\lambda_3}b_{\lambda_2-\lambda_3}
      -b_{-\lambda_3}b_{\lambda_3}.
  \end{equation}
  If $b\neq0$, set $X=e_{\alpha}/(3b)$, $Y=3be_{-\alpha}$, and $x=[X,y]$.
  Then $(X,x,h,y,Y)$ is a minimal $\mathfrak{osp}(1\vert 2)$-quintuple.
\end{prop}

\begin{proof}
  By the form of $y$, $[h,y]=-y$.  The relations in \Cref{sec:commutation_relations_for_F(4)} give
  \[-\frac{1}{2}[y,y]=3be_{-\alpha}.\]
  If $b\neq0$, the prescribed $(X,h,Y)$ is a minimal $\mathfrak{sl}_2$-triple by \Cref{prop:minimal_triple_entries} and \Cref{rem:scalar_multiple_minimal_triple}.
  These are the conditions in \Cref{lem:minimality_criterion}, so the assertion follows.
\end{proof}

In particular, the extensions of $(e_{\alpha}/3,h,3e_{-\alpha})$ are precisely the quintuples in \Cref{prop:case_1_of_F(4)} with $b=1$.

\begin{prop}\label{prop:F4_case_1_conjugacy}
  Suppose $(e_{\alpha}/3,x_1,h,y_1,3e_{-\alpha})$ and $(e_{\alpha}/3,x_2,h,y_2,3e_{-\alpha})$ are minimal $\mathfrak{osp}(1\vert 2)$-quintuples whose entries $y_1$ and $y_2$ are as in \Cref{prop:case_1_of_F(4)} with $b=1$.
  Then the two quintuples are conjugate.
\end{prop}

\begin{proof}
  Order the basis vectors in \Cref{prop:case_1_of_F(4)} as follows:
  \[\begin{aligned}
    &w_0=v_{-\lambda_3}^{-},\quad w_1=v_{-\lambda_2+\lambda_3}^{-},\quad w_2=v_{-\lambda_1+\lambda_2-\lambda_3}^{-},\quad w_3=v_{\lambda_1-\lambda_3}^{-},\\
    &w_4=v_{-\lambda_1+\lambda_3}^{-},\quad w_5=v_{\lambda_1-\lambda_2+\lambda_3}^{-},\quad w_6=v_{\lambda_2-\lambda_3}^{-},\quad w_7=v_{\lambda_3}^{-}.
  \end{aligned}\]
  Write $y=\sum_{i=0}^7b_iw_i$.
  The condition $b=1$ is
  \begin{equation}\label{eq:F4_case_1_quadratic}
    b_3b_4-b_2b_5+b_1b_6-b_0b_7=1.
  \end{equation}
  Set
  \[(\gamma_1,\gamma_2,\gamma_3,\gamma_4,\gamma_5,\gamma_6)=(\alpha_3,\alpha_5,\alpha_6,\alpha_7,\alpha_8,\alpha_9).\]
  The relations in \Cref{sec:commutation_relations_for_F(4)} give
  \[[e_{\gamma_i},w_0]=w_i\qquad\text{and}\qquad[e_{-\gamma_i},w_i]=w_0\qquad(1\leq i\leq6).\]
  Moreover, for $1\leq i\leq6$ and $1\leq j\leq7$, we have
  \[[e_{\gamma_i},w_j]\in\operatorname{span}_{\mathbb{C}}\{w_{i+1},\dotsc,w_7\}.\]

  Suppose first that $b_0=0$.
  Equation \eqref{eq:F4_case_1_quadratic} implies that $b_i\neq0$ for some $1\leq i\leq6$.
  For such an $i$, $[e_{-\gamma_i},w_i]=w_0$, and no other $w_j$ contributes to the $w_0$-coefficient under a single application of $\operatorname{ad}(e_{-\gamma_i})$.
  Thus the $w_0$-coefficient of $\sigma_{-\gamma_i,c}(y)$ is a polynomial in $c$ with nonzero linear term $b_i c$.
  Choosing $c$ away from its zeros, we may assume $b_0\neq0$.
  Successively for $i=1,\dotsc,6$, apply $\sigma_{\gamma_i,c_i}$.
  At the $i$-th step, this automorphism changes the current coefficient of $w_i$ by $c_i b_0$ and leaves the coefficients of $w_0,\dotsc,w_{i-1}$ unchanged.
  We can therefore choose $c_i$ to cancel the coefficient of $w_i$ without reintroducing any coefficient canceled at an earlier step.
  We thereby obtain an element of the form
  \[\widetilde{y}=b_0w_0+\widetilde{b}_7w_7.\]
  All the automorphisms used fix $(e_{\alpha}/3,h,3e_{-\alpha})$ by \Cref{lem:root_exponentials_fix_triple}; hence \eqref{eq:F4_case_1_quadratic} remains valid and gives $\widetilde{b}_7=-b_0^{-1}$.
  Finally, since $[h_3,w_0]=-w_0$ and $[h_3,w_7]=w_7$, the diagonal automorphism $\tau_{\alpha_3,b_0}$ sends $\widetilde{y}$ to $w_0-w_7$.
  Thus both $y_1$ and $y_2$ are conjugate to $w_0-w_7$ by automorphisms fixing $(e_{\alpha}/3,h,3e_{-\alpha})$.
  It follows that the given two quintuples are conjugate.
\end{proof}

For a quintuple extending the second triple, its fourth entry $y$ satisfies $[h_9,y]=-y$; the weights in \Cref{fig:F4-weight-diagram} then force $y$ to lie in the span of $v_{-\lambda_3}^{\pm}$ and $v_{-\lambda_2+\lambda_3}^{\pm}$.

\begin{prop}\label{prop:case_2_of_F(4)}
  Let $a_1,a_2,a_3,a_4\in\mathbb{C}$, and set
  \begin{equation}\label{eq:F4_case_2_element}
    y=a_1v_{-\lambda_3}^{-}+a_2v_{-\lambda_3}^{+}
      +a_3v_{-\lambda_2+\lambda_3}^{-}+a_4v_{-\lambda_2+\lambda_3}^{+}
    \qquad\text{and}\qquad
    a:=a_1a_4-a_2a_3.
  \end{equation}
  If $a\neq0$, set $X=e_{\alpha_9}/(2a)$, $Y=2ae_{-\alpha_9}$, and $x=[X,y]$.
  Then $(X,x,h_9,y,Y)$ is a minimal $\mathfrak{osp}(1\vert 2)$-quintuple.
\end{prop}

\begin{proof}
  By the form of $y$, $[h_9,y]=-y$.  The relations in \Cref{sec:commutation_relations_for_F(4)} give
  \[-\frac{1}{2}[y,y]=2ae_{-\alpha_9}.\]
  If $a\neq0$, the prescribed $(X,h_9,Y)$ is a minimal $\mathfrak{sl}_2$-triple by \Cref{prop:minimal_triple_entries} and \Cref{rem:scalar_multiple_minimal_triple}.
  These are the conditions in \Cref{lem:minimality_criterion}, so the assertion follows.
\end{proof}

In particular, the extensions of $(e_{\alpha_9}/2,h_9,2e_{-\alpha_9})$ are precisely the quintuples in \Cref{prop:case_2_of_F(4)} with $a=1$.

\begin{prop}\label{prop:F4_case_2_conjugacy}
  Suppose $(e_{\alpha_9}/2,x_1,h_9,y_1,2e_{-\alpha_9})$ and $(e_{\alpha_9}/2,x_2,h_9,y_2,2e_{-\alpha_9})$ are minimal $\mathfrak{osp}(1\vert 2)$-quintuples whose entries $y_1$ and $y_2$ are as in \eqref{eq:F4_case_2_element} with $a=1$.
  Then the two quintuples are conjugate.
\end{prop}

\begin{proof}
  Identify the element in \eqref{eq:F4_case_2_element} with the matrix
  \[A=\begin{pmatrix}a_1&a_2\\a_3&a_4\end{pmatrix}.\]
  The normalization $a=1$ means that $\det A=1$.
  Since $\alpha_3(h_9)=0$, \Cref{lem:root_exponentials_fix_triple} shows that $\sigma_{\alpha_3,c}$ and $\sigma_{-\alpha_3,c}$ fix the triple $(e_{\alpha_9}/2,h_9,2e_{-\alpha_9})$ pointwise.
  From the relations in \Cref{sec:commutation_relations_for_F(4)}, we have
  \[\sigma_{\alpha_3,c}\cdot A
    =\begin{pmatrix}1&0\\c&1\end{pmatrix}A,
    \qquad
    \sigma_{-\alpha_3,c}\cdot A
    =\begin{pmatrix}1&c\\0&1\end{pmatrix}A.\]
  These are left multiplications by the two types of elementary unipotent matrices, which generate $\operatorname{SL}_2(\mathbb{C})$.
  In particular, if $A_1,A_2\in\operatorname{SL}_2(\mathbb{C})$ correspond to $y_1,y_2$, then $A_2A_1^{-1}$ is a product of these elementary matrices, so $y_1$ and $y_2$ are conjugate by automorphisms fixing the triple.
  The same automorphisms therefore conjugate the quintuples.
\end{proof}

\subsection{The Centralizers and Complementary Submodules in \texorpdfstring{$F(4)$}{F(4)}}

For each conjugacy class found in the preceding subsection, we choose a representative $\mathcal{Q}$ with even triple $\mathcal{T}$ and compute the centralizers $\mathfrak{g}^{\mathcal{T}}$ and $\mathfrak{g}^{\mathcal{Q}}$.
We also construct a complementary submodule $M_{\mathcal{Q}}^{\mathfrak{g}}$ and determine its isomorphism type.

\begin{prop}\label{prop:F4_first_centralizers}
  Consider the minimal $\mathfrak{osp}(1\vert 2)$-quintuple $\mathcal{Q}$ given by
  \[\mathcal{Q}=\left(\frac{1}{3}e_{\alpha},
      \frac{1}{3}\bigl(v_{-\lambda_3}^{+}-v_{\lambda_3}^{+}\bigr),
      h,
      v_{-\lambda_3}^{-}-v_{\lambda_3}^{-},
      3e_{-\alpha}\right),\]
  and write $\mathcal{T}$ for its even triple.
  Then
  \[\mathfrak{g}^{\mathcal{T}}=\mathfrak{a}_2\cong\mathfrak{so}_7,\qquad
    \mathfrak{g}^{\mathcal{Q}}\cong G_2.\]
  Moreover, there exists a complementary submodule $M_{\mathcal{Q}}^{\mathfrak{g}}$, and
  \[M_{\mathcal{Q}}^{\mathfrak{g}}\cong V_{G_2}(\varpi_1),\]
  where $\varpi_1$ is the fundamental weight associated with the short simple root of $G_2$.
\end{prop}

\begin{proof}
  The entries of $\mathcal{T}$ span $\mathfrak{a}_1$, so their joint centralizer in $\mathfrak{a}_1$ is zero.
  Their joint centralizer in $\mathfrak{g}_{\bar{1}}\cong V_1\otimes V_2$ is also zero, because $V_1$ is the natural $\mathfrak{sl}_2$-module.
  Since $\mathfrak{a}_1$ commutes with $\mathfrak{a}_2$, we obtain
  \[\mathfrak{g}^{\mathcal{T}}=\mathfrak{a}_2.\]

  To compute $\mathfrak{g}^{\mathcal{Q}}$, write $y$ for the fourth entry of $\mathcal{Q}$ and let \[g=\sum_{i=1}^{9}a_ie_{\alpha_i}+\sum_{i=1}^{9}b_ie_{-\alpha_i}+\sum_{i=1}^{3}c_ih_i\] be an arbitrary element of $\mathfrak{g}^{\mathcal{T}}$, where $a_i,b_i,c_i\in\mathbb{C}$.
  Then we have
  \begin{multline*}
  [g,y]=-c_3v_{-\lambda_3}^{-}+(a_3-b_9)v_{-\lambda_2+\lambda_3}^{-}+(a_5+b_8)v_{-\lambda_1+\lambda_2-\lambda_3}^{-}+(a_6-b_7)v_{\lambda_1-\lambda_3}^{-}\\
  +(a_7-b_6)v_{-\lambda_1+\lambda_3}^{-}+(a_8+b_5)v_{\lambda_1-\lambda_2+\lambda_3}^{-}+(a_9-b_3)v_{\lambda_2-\lambda_3}^{-}-c_3v_{\lambda_3}^{-}.
  \end{multline*}
  Thus, $\mathfrak{g}^{\mathcal{Q}}=\{g\in\mathfrak{g}^{\mathcal{T}}\mid[g,y]=0\}$ is spanned by
  \[\begin{aligned}
    &e_{\alpha_5}-e_{-\alpha_8}, &&\qquad e_{\alpha_1}, &&\qquad e_{\alpha_6}+e_{-\alpha_7}, &&\qquad e_{\alpha_9}+e_{-\alpha_3}, &&\qquad e_{\alpha_2}, &&\qquad -e_{\alpha_4},\\
    &e_{-\alpha_5}-e_{\alpha_8}, &&\qquad e_{-\alpha_1}, &&\qquad e_{-\alpha_6}+e_{\alpha_7}, &&\qquad e_{-\alpha_9}+e_{\alpha_3}, &&\qquad e_{-\alpha_2}, &&\qquad -e_{-\alpha_4},\\
    &-h_1+h_2, &&\qquad h_1.
  \end{aligned}\]
  Denote the six elements in the first row by $u_1,\dotsc,u_6$, the six elements in the second row by $u_{-1},\dotsc,u_{-6}$, and set $k_1=-h_1+h_2$ and $k_2=h_1$.
  The relations in \Cref{sec:commutation_relations_for_F(4)} give
  \[[u_1,u_{-1}]=k_1,\qquad [u_2,u_{-2}]=k_2,\qquad
    \bigl([k_i,u_j]\bigr)_{1\leq i,j\leq2}
    =\begin{pmatrix}2u_1&-3u_2\\-u_1&2u_2\end{pmatrix}.\]
  Moreover,
  \[[u_2,u_1]=u_3,\qquad [u_1,u_3]=2u_4,\qquad
    [u_1,u_4]=3u_5,\qquad [u_2,u_5]=u_6.\]
  The coefficient matrix $\begin{pmatrix}2&-3\\-1&2\end{pmatrix}$ is the Cartan matrix of type $G_2$, with the short simple root listed first.
  The four subsequent brackets produce the remaining positive root vectors; the negative root vectors $u_{-3},\dotsc,u_{-6}$ satisfy the corresponding relations.
  The fourteen listed vectors are linearly independent and exhaust $\mathfrak{g}^{\mathcal{Q}}$, so these brackets identify its root space decomposition and give $\mathfrak{g}^{\mathcal{Q}}\cong G_2$.

  Now let
  \[\begin{aligned}
    M_{\mathcal{Q}}^{\mathfrak{g}}&=\operatorname{span}_{\mathbb{C}}\{-e_{\alpha_3}+2e_{-\alpha_9},\; 2e_{\alpha_7}-e_{-\alpha_6},\; 2e_{\alpha_8}+e_{-\alpha_5},\\
    &\qquad\qquad\quad 2h_1+4h_2+3h_3,\; -2e_{-\alpha_8}-e_{\alpha_5},\; 2e_{-\alpha_7}-e_{\alpha_6},\; e_{-\alpha_3}-2e_{\alpha_9}\}.
  \end{aligned}\]
  The seven vectors spanning $M_{\mathcal{Q}}^{\mathfrak{g}}$ lie in $\mathfrak{g}^{\mathcal{T}}$ and are linearly independent modulo $\mathfrak{g}^{\mathcal{Q}}$; denote them, in order, by $m_1,\dotsc,m_7$.
  Since $\dim\mathfrak{g}^{\mathcal{T}}=21$, $\dim\mathfrak{g}^{\mathcal{Q}}=14$, and $\dim M_{\mathcal{Q}}^{\mathfrak{g}}=7$, this independence gives $\mathfrak{g}^{\mathcal{T}}=\mathfrak{g}^{\mathcal{Q}}\oplus M_{\mathcal{Q}}^{\mathfrak{g}}$ as vector spaces.
  The relations in \Cref{sec:commutation_relations_for_F(4)} give
  \[[u_{\pm i},M_{\mathcal{Q}}^{\mathfrak{g}}]\subseteq M_{\mathcal{Q}}^{\mathfrak{g}}\qquad(i=1,2).\]
  Since $u_{\pm1}$ and $u_{\pm2}$ generate $\mathfrak{g}^{\mathcal{Q}}$, this proves that $M_{\mathcal{Q}}^{\mathfrak{g}}$ is a $\mathfrak{g}^{\mathcal{Q}}$-submodule.
  The actions on the seven basis vectors include the following:
  \[\begin{gathered}
    {}[u_1,m_1]=m_2,\qquad [u_2,m_2]=[u_3,m_1]=m_3,\\
    [u_1,m_3]=[u_4,m_1]=m_4,\qquad
    [u_1,m_4]=2m_5,\qquad [u_5,m_1]=m_5,\\
    [u_2,m_5]=[u_6,m_1]=m_6,\qquad
    [u_1,m_6]=[u_6,m_2]=m_7.
  \end{gathered}\]
  Put $\mathfrak{t}=\operatorname{span}_{\mathbb{C}}\{k_1,k_2\}$, a Cartan subalgebra of $\mathfrak{g}^{\mathcal{Q}}$.
  Let $\varpi_1,\varpi_2\in\mathfrak{t}^*$ be the fundamental weights dual to the coroots $k_1,k_2$, respectively; thus $\varpi_i(k_j)=\delta_{ij}$ for $i,j=1,2$.
  The $\mathfrak{t}$-weights of $m_1,\dotsc,m_7$, in order, are
  \[-\varpi_1,\quad \varpi_1-\varpi_2,\quad -2\varpi_1+\varpi_2,\quad 0,\quad
    2\varpi_1-\varpi_2,\quad -\varpi_1+\varpi_2,\quad \varpi_1.\]
  Since $[u_1,m_7]=[u_2,m_7]=0$, the vector $m_7$ has highest weight $\varpi_1$.
  Since finite-dimensional $\mathfrak{g}^{\mathcal{Q}}$-modules are completely reducible and $\dim V_{\mathfrak{g}^{\mathcal{Q}}}(\varpi_1)=7=\dim M_{\mathcal{Q}}^{\mathfrak{g}}$, it follows that $M_{\mathcal{Q}}^{\mathfrak{g}}\cong V_{\mathfrak{g}^{\mathcal{Q}}}(\varpi_1)$.
\end{proof}

\begin{prop}\label{prop:F4_second_centralizers}
  Consider the minimal $\mathfrak{osp}(1\vert 2)$-quintuple $\mathcal{Q}$ given by
  \[\mathcal{Q}=\left(\frac{1}{2}e_{\alpha_9},
      \frac{1}{2}\bigl(v_{\lambda_2-\lambda_3}^{+}-v_{\lambda_3}^{-}\bigr),
      h_9,
      v_{-\lambda_3}^{+}-v_{-\lambda_2+\lambda_3}^{-},
      2e_{-\alpha_9}\right),\]
  and write $\mathcal{T}$ for its even triple.
  Then
  \[\mathfrak{g}^{\mathcal{T}}\cong D(2,1;2),\qquad
    \mathfrak{g}^{\mathcal{Q}}\cong\mathfrak{sl}_2\oplus\mathfrak{osp}(1\vert 2).\]
  Moreover, there exists a complementary submodule $M_{\mathcal{Q}}^{\mathfrak{g}}$, and
  \[M_{\mathcal{Q}}^{\mathfrak{g}}\cong\mathbb{C}^3\otimes\mathbb{C}^{1\vert 2}.\]
  Here, $\mathbb{C}^3$ is the $3$-dimensional irreducible module for the $\mathfrak{sl}_2$-summand of $\mathfrak{g}^{\mathcal{Q}}$, and $\mathbb{C}^{1\vert 2}$ is the irreducible $\mathfrak{osp}(1\vert 2)$-module of dimension $(1\vert 2)$.
\end{prop}

\begin{proof}
  The bracket relations with $e_{\pm\alpha_9}$ give
  \[\begin{aligned}
    \mathfrak{g}^{\mathcal{T}}&=\operatorname{span}_{\mathbb{C}}\left\{e_{\alpha},\,h,\,e_{-\alpha};\,e_{\alpha_1},\,h_1,\,e_{-\alpha_1};\,e_{\alpha_3},\,h_3,\,e_{-\alpha_3}\right\}\\
    &\qquad\oplus\operatorname{span}_{\mathbb{C}}\left\{v_{\lambda_1-\lambda_2+\lambda_3}^{\pm},\,v_{-\lambda_1+\lambda_3}^{\pm},\,v_{\lambda_1-\lambda_3}^{\pm},\,v_{-\lambda_1+\lambda_2-\lambda_3}^{\pm}\right\}.
  \end{aligned}\]
  Put
  \[\mathfrak{s}_0=\operatorname{span}_{\mathbb{C}}\{e_{\alpha},h,e_{-\alpha}\},\qquad
    \mathfrak{s}_1=\operatorname{span}_{\mathbb{C}}\{e_{\alpha_1},h_1,e_{-\alpha_1}\},\qquad
    \mathfrak{s}_3=\operatorname{span}_{\mathbb{C}}\{e_{\alpha_3},h_3,e_{-\alpha_3}\}.\]
  Then $(\mathfrak{g}^{\mathcal{T}})_{\bar{0}}=\mathfrak{s}_0\oplus\mathfrak{s}_1\oplus\mathfrak{s}_3$, where the three summands are commuting copies of $\mathfrak{sl}_2$.
  Let $U_i$ be the natural $2$-dimensional module for $\mathfrak{s}_i$, for $i=0,1,3$.
  Relative to $(h,h_1,h_3)$, the eight odd basis vectors have the distinct weights $(\pm1,\pm1,\pm1)$, and the relations in \Cref{sec:commutation_relations_for_F(4)} identify $(\mathfrak{g}^{\mathcal{T}})_{\bar{1}}$ with $U_0\otimes U_1\otimes U_3$ as a module for $(\mathfrak{g}^{\mathcal{T}})_{\bar{0}}$.
  The relations
  \[\begin{aligned}
    \left[v_{\lambda_1-\lambda_2+\lambda_3}^{+},v_{-\lambda_1+\lambda_2-\lambda_3}^{+}\right]&=-3e_{\alpha},\\
    \left[v_{\lambda_1-\lambda_2+\lambda_3}^{+},v_{\lambda_1-\lambda_3}^{-}\right]&=2e_{\alpha_1},\\
    \left[v_{\lambda_1-\lambda_3}^{+},v_{-\lambda_1+\lambda_2-\lambda_3}^{-}\right]&=-e_{-\alpha_3}
  \end{aligned}\]
  give nonzero root vectors in the three simple ideals of $(\mathfrak{g}^{\mathcal{T}})_{\bar{0}}$, respectively.
  By the super Jacobi identity, the span of odd-odd brackets is an ideal of $(\mathfrak{g}^{\mathcal{T}})_{\bar{0}}$; hence
  \[\bigl[(\mathfrak{g}^{\mathcal{T}})_{\bar{1}},(\mathfrak{g}^{\mathcal{T}})_{\bar{1}}\bigr]=(\mathfrak{g}^{\mathcal{T}})_{\bar{0}}.\]
  Since $(\mathfrak{g}^{\mathcal{T}})_{\bar{1}}$ is a faithful irreducible $(\mathfrak{g}^{\mathcal{T}})_{\bar{0}}$-module, \Cref{lem:simplicity_from_odd_part} shows that $\mathfrak{g}^{\mathcal{T}}$ is simple.
  Together with the descriptions of its even and odd parts, the classification of finite-dimensional simple Lie superalgebras now gives $\mathfrak{g}^{\mathcal{T}}\cong D(2,1;\gamma)$ for some $\gamma$.

  It remains to determine $\gamma$.
  For each $i\in\{0,1,3\}$, choose a weight basis $\{u_i^+,u_i^-\}$ of $U_i$ such that the positive root vector of $\mathfrak{s}_i$ sends $u_i^-$ to $u_i^+$, and normalize the invariant alternating form $\langle\,\cdot\,,\,\cdot\,\rangle_i$ on $U_i$ by $\langle u_i^+,u_i^-\rangle_i=1$.
  Match each odd basis vector with the tensor basis vector of the same weight under the identification above.
  For $u,v\in U_i$, define $P_i(u,v)\in\mathfrak{s}_i\cong\mathfrak{sp}(U_i)$ by
  \[P_i(u,v)w=\langle v,w\rangle_i u+\langle u,w\rangle_i v\qquad(w\in U_i).\]
  In the standard tensor-product realization of $D(2,1;\gamma)$, let $c_i$ be the coefficient of the $\mathfrak{s}_i$-component of the odd-odd bracket obtained by applying $P_i$ to the $U_i$-factors and contracting the other two pairs with their alternating forms.
  Up to a common nonzero scalar and a permutation of the three $\mathfrak{sl}_2$-factors, the coefficient triple is $(1,\gamma,-1-\gamma)$.
  The three relations above isolate the $\mathfrak{s}_0$-, $\mathfrak{s}_1$-, and $\mathfrak{s}_3$-components, respectively.
  Since $P_i(u_i^+,u_i^+)$ is twice the positive root vector of $\mathfrak{s}_i$, whereas $P_i(u_i^-,u_i^-)$ is minus twice its negative root vector, they give
  \[2c_0=-3,\qquad 2c_1=2,\qquad -2c_3=-1,\qquad
    (c_0,c_1,c_3)=\frac{1}{2}(-3,2,1).\]
  After permuting the first and third factors, this is $\frac{1}{2}(1,2,-3)=\frac{1}{2}(1,2,-1-2)$.
  It follows that $\mathfrak{g}^{\mathcal{T}}\cong D(2,1;2)$.

  To compute $\mathfrak{g}^{\mathcal{Q}}$, write $y$ for the fourth entry of $\mathcal{Q}$.
  For an arbitrary $g_{\bar{0}}\in(\mathfrak{g}^{\mathcal{T}})_{\bar{0}}$, write
  \[g_{\bar{0}}=c_1e_{\alpha}+c_2e_{-\alpha}+c_3h+c_4e_{\alpha_1}+c_5e_{-\alpha_1}+c_6h_1+c_7e_{\alpha_3}+c_8e_{-\alpha_3}+c_9h_3,\]
  with $c_1,\dotsc,c_9\in\mathbb{C}$.
  Then
  \[[g_{\bar{0}},y]=(c_3-c_9)v_{-\lambda_3}^{+}+(c_2-c_8)v_{-\lambda_3}^{-}-(c_1-c_7)v_{-\lambda_2+\lambda_3}^{+}+(c_3-c_9)v_{-\lambda_2+\lambda_3}^{-}.\]
  Thus $[g_{\bar{0}},y]=0$ precisely when $c_1=c_7$, $c_2=c_8$, and $c_3=c_9$.
  Hence
  \[\begin{aligned}
    (\mathfrak{g}^{\mathcal{Q}})_{\bar{0}}&=\left\{g_{\bar{0}}\in(\mathfrak{g}^{\mathcal{T}})_{\bar{0}}\mid[g_{\bar{0}},y]=0\right\}\\
    &=\operatorname{span}_{\mathbb{C}}\left\{e_{\alpha}+e_{\alpha_3},\,e_{-\alpha}+e_{-\alpha_3},\,h+h_3;\,e_{\alpha_1},\,e_{-\alpha_1},\,h_1\right\}\\
    &\cong\mathfrak{sl}_2\oplus\mathfrak{sl}_2.
  \end{aligned}\]
  Similarly, for an arbitrary $g_{\bar{1}}\in(\mathfrak{g}^{\mathcal{T}})_{\bar{1}}$, write
  \begin{multline*}
    g_{\bar{1}}=d_1v_{-\lambda_1+\lambda_2-\lambda_3}^{+}+d_2v_{-\lambda_1+\lambda_2-\lambda_3}^{-}+d_3v_{\lambda_1-\lambda_3}^{+}+d_4v_{\lambda_1-\lambda_3}^{-}\\
    +d_5v_{-\lambda_1+\lambda_3}^{+}+d_6v_{-\lambda_1+\lambda_3}^{-}+d_7v_{\lambda_1-\lambda_2+\lambda_3}^{+}+d_8v_{\lambda_1-\lambda_2+\lambda_3}^{-},
  \end{multline*}
  with $d_1,\dotsc,d_8\in\mathbb{C}$.
  Then
  \[[g_{\bar{1}},y]=-2d_7e_{-\alpha_2}-2d_5e_{-\alpha_4}-(d_3+d_8)e_{-\alpha_5}+(d_1+d_6)e_{-\alpha_6}+2d_4e_{-\alpha_7}-2d_2e_{-\alpha_8}.\]
  The vanishing condition is $d_2=d_4=d_5=d_7=0$, $d_6=-d_1$, and $d_8=-d_3$.
  Hence
  \[\begin{aligned}
    (\mathfrak{g}^{\mathcal{Q}})_{\bar{1}}&=\left\{g_{\bar{1}}\in(\mathfrak{g}^{\mathcal{T}})_{\bar{1}}\mid[g_{\bar{1}},y]=0\right\}\\
    &=\operatorname{span}_{\mathbb{C}}\left\{\frac{1}{\sqrt{2}}\left(v_{-\lambda_1+\lambda_2-\lambda_3}^{+}-v_{-\lambda_1+\lambda_3}^{-}\right),\,\frac{1}{\sqrt{2}}\left(v_{\lambda_1-\lambda_3}^{+}-v_{\lambda_1-\lambda_2+\lambda_3}^{-}\right)\right\}.
  \end{aligned}\]
  Set
  \[\begin{gathered}
    X'=e_{\alpha}+e_{\alpha_3},\qquad H'=h+h_3,\qquad Y'=e_{-\alpha}+e_{-\alpha_3},\\
    X''=e_{\alpha_1},\qquad H''=h_1,\qquad Y''=e_{-\alpha_1},\\
    x'':=\frac{1}{\sqrt{2}}\left(v_{\lambda_1-\lambda_3}^{+}-v_{\lambda_1-\lambda_2+\lambda_3}^{-}\right),\qquad
    y'':=\frac{1}{\sqrt{2}}\left(v_{-\lambda_1+\lambda_2-\lambda_3}^{+}-v_{-\lambda_1+\lambda_3}^{-}\right).
  \end{gathered}\]
  The relations in \Cref{sec:commutation_relations_for_F(4)} show that $(X',H',Y')$ is an $\mathfrak{sl}_2$-triple, that $(X'',x'',H'',y'',Y'')$ is an $\mathfrak{osp}(1\vert 2)$-quintuple, and that their spans commute.
  Since these elements together form a basis of $\mathfrak{g}^{\mathcal{Q}}$, it follows that
  \[\mathfrak{g}^{\mathcal{Q}}\cong\mathfrak{sl}_2\oplus\mathfrak{osp}(1\vert 2).\]

  Define
  \[\begin{aligned}
    M_{\mathcal{Q}}^{\mathfrak{g}}&=\operatorname{span}_{\mathbb{C}}\left\{v_{-\lambda_1+\lambda_2-\lambda_3}^{-},\,v_{-\lambda_1+\lambda_2-\lambda_3}^{+}+v_{-\lambda_1+\lambda_3}^{-},\,v_{-\lambda_1+\lambda_3}^{+}\right\}\\
    &\quad\oplus\operatorname{span}_{\mathbb{C}}\left\{-\frac{1}{\sqrt{2}}(3e_{-\alpha}+e_{-\alpha_3}),\,-\frac{1}{\sqrt{2}}(3h+h_3),\,\frac{1}{\sqrt{2}}(3e_{\alpha}+e_{\alpha_3})\right\}\\
    &\quad\oplus\operatorname{span}_{\mathbb{C}}\left\{v_{\lambda_1-\lambda_3}^{-},\,v_{\lambda_1-\lambda_3}^{+}+v_{\lambda_1-\lambda_2+\lambda_3}^{-},\,v_{\lambda_1-\lambda_2+\lambda_3}^{+}\right\}.
  \end{aligned}\]
  The nine vectors spanning $M_{\mathcal{Q}}^{\mathfrak{g}}$ lie in $\mathfrak{g}^{\mathcal{T}}$ and are linearly independent modulo $\mathfrak{g}^{\mathcal{Q}}$.
  Since $\dim\mathfrak{g}^{\mathcal{T}}=17$, $\dim\mathfrak{g}^{\mathcal{Q}}=8$, and $\dim M_{\mathcal{Q}}^{\mathfrak{g}}=9$, we obtain $\mathfrak{g}^{\mathcal{T}}=\mathfrak{g}^{\mathcal{Q}}\oplus M_{\mathcal{Q}}^{\mathfrak{g}}$ as vector spaces.
  Regard the three spanning sets as rows indexed, from top to bottom, by $b=-1,0,1$, and the vectors in each row as columns indexed, from left to right, by $a=-2,0,2$.
  Write $m_{a,b}$ for the vector in position $(a,b)$.
  Direct calculation gives
  \[[H',m_{a,b}]=a m_{a,b},\qquad [H'',m_{a,b}]=b m_{a,b}.\]
  Furthermore, for every $b\in\{-1,0,1\}$,
  \[\begin{aligned}
    {}[X',m_{-2,b}]&=m_{0,b}, & [X',m_{0,b}]&=2m_{2,b},\\
    [Y',m_{2,b}]&=m_{0,b}, & [Y',m_{0,b}]&=2m_{-2,b},
  \end{aligned}\]
  while, for every $a\in\{-2,0,2\}$,
  \[\begin{aligned}
    {}[x'',m_{a,-1}]&=m_{a,0}, & [x'',m_{a,0}]&=m_{a,1},\\
    [y'',m_{a,1}]&=-m_{a,0}, & [y'',m_{a,0}]&=m_{a,-1}.
  \end{aligned}\]
  All remaining actions of $X',Y',x'',y''$ on the vectors $m_{a,b}$ are zero.
  Since these four elements generate $\mathfrak{g}^{\mathcal{Q}}$, the relations show that $M_{\mathcal{Q}}^{\mathfrak{g}}$ is a $\mathfrak{g}^{\mathcal{Q}}$-submodule isomorphic to $\mathbb{C}^{3}\otimes\mathbb{C}^{1\vert 2}$.
\end{proof}

\begin{rem}
  The calculations in \Cref{prop:F4_first_centralizers,prop:F4_second_centralizers} show that, in each case, $\mathfrak{g}^{\mathcal{T}}$ is a direct sum of pairwise nonisomorphic irreducible $\mathfrak{g}^{\mathcal{Q}}$-submodules.
  Consequently, the complementary submodule $M_{\mathcal{Q}}^{\mathfrak{g}}$ constructed in each case is in fact the unique $\mathfrak{g}^{\mathcal{Q}}$-submodule of $\mathfrak{g}^{\mathcal{T}}$ complementary to $\mathfrak{g}^{\mathcal{Q}}$.
\end{rem}

\subsection{Summary for \texorpdfstring{$F(4)$}{F(4)}}

We summarize the results of this section.

\begin{thm}\label{thm:main_F(4)}
  Let $\mathfrak{g}=F(4)$, and let $\mathfrak{a}_1$ and $\mathfrak{a}_2$ be the simple ideals of $\mathfrak{g}_{\bar{0}}$ isomorphic to $\mathfrak{sl}_2$ and $\mathfrak{so}_7$, respectively.
  There are exactly two conjugacy classes of minimal $\mathfrak{osp}(1\vert 2)$-quintuples in $\mathfrak{g}$.
  Representatives $\mathcal{Q}_i=(X_i,x_i,H_i,y_i,Y_i)$, $i=1,2$, of these classes are given by
  \begin{align*}
    &X_1=\frac{1}{3}e_{\alpha},
    &&x_1=\frac{1}{3}\bigl(v_{-\lambda_3}^{+}-v_{\lambda_3}^{+}\bigr),
    &&H_1=h,\\
    &y_1=v_{-\lambda_3}^{-}-v_{\lambda_3}^{-},
    &&Y_1=3e_{-\alpha},\\[2pt]
    &X_2=\frac{1}{2}e_{\alpha_9},
    &&x_2=\frac{1}{2}\bigl(v_{\lambda_2-\lambda_3}^{+}-v_{\lambda_3}^{-}\bigr),
    &&H_2=h_9,\\
    &y_2=v_{-\lambda_3}^{+}-v_{-\lambda_2+\lambda_3}^{-},
    &&Y_2=2e_{-\alpha_9}.
  \end{align*}
  Let $\mathcal{T}_i=(X_i,H_i,Y_i)$ be the even triple of $\mathcal{Q}_i$ for $i=1,2$.
  If $\varpi_1$ is the fundamental weight associated with the short simple root of $G_2$, then the associated centralizers and complementary submodules are
  \begin{align*}
    &\mathfrak{g}^{\mathcal{T}_1}=\mathfrak{a}_2\cong\mathfrak{so}_7,
    &&\mathfrak{g}^{\mathcal{Q}_1}\cong G_2,
    &&M_{\mathcal{Q}_1}^{\mathfrak{g}}\cong V_{G_2}(\varpi_1),\\
    &\mathfrak{g}^{\mathcal{T}_2}\cong D(2,1;2),
    &&\mathfrak{g}^{\mathcal{Q}_2}\cong\mathfrak{sl}_2\oplus\mathfrak{osp}(1\vert 2),
    &&M_{\mathcal{Q}_2}^{\mathfrak{g}}\cong\mathbb{C}^3\otimes\mathbb{C}^{1\vert 2}.
  \end{align*}
  In the last isomorphism, $\mathbb{C}^3$ and $\mathbb{C}^{1\vert 2}$ denote the irreducible modules for the respective summands $\mathfrak{sl}_2$ and $\mathfrak{osp}(1\vert 2)$ of $\mathfrak{g}^{\mathcal{Q}_2}$.
\end{thm}

\begin{proof}
  The representatives and their centralizers and complementary submodules are given by \Cref{prop:F4_first_centralizers,prop:F4_second_centralizers}.
  The classification of minimal $\mathfrak{sl}_2$-triples in \Cref{subsec:minimal_nilpotent_orbits}, together with the normalization at the start of \Cref{subsec:F4_conjugacy}, shows that every minimal $\mathfrak{osp}(1\vert 2)$-quintuple is conjugate to an extension of $\mathcal{T}_1$ or $\mathcal{T}_2$.
  The extensions of each triple are conjugate by \Cref{prop:F4_case_1_conjugacy,prop:F4_case_2_conjugacy}; hence these representatives give exactly the two conjugacy classes.
\end{proof}

\section{Minimal \texorpdfstring{$\mathfrak{osp}(1\vert 2)$}{osp(1|2)}-quintuples and Their Centralizers in \texorpdfstring{$G(3)$}{G(3)}}\label{sec:G3}

Throughout this section, let $\mathfrak{g}=G(3)$.

\subsection{An Explicit Basis for \texorpdfstring{$G(3)$}{G(3)}}\label{subsec:explicit_basis_for_G(3)}

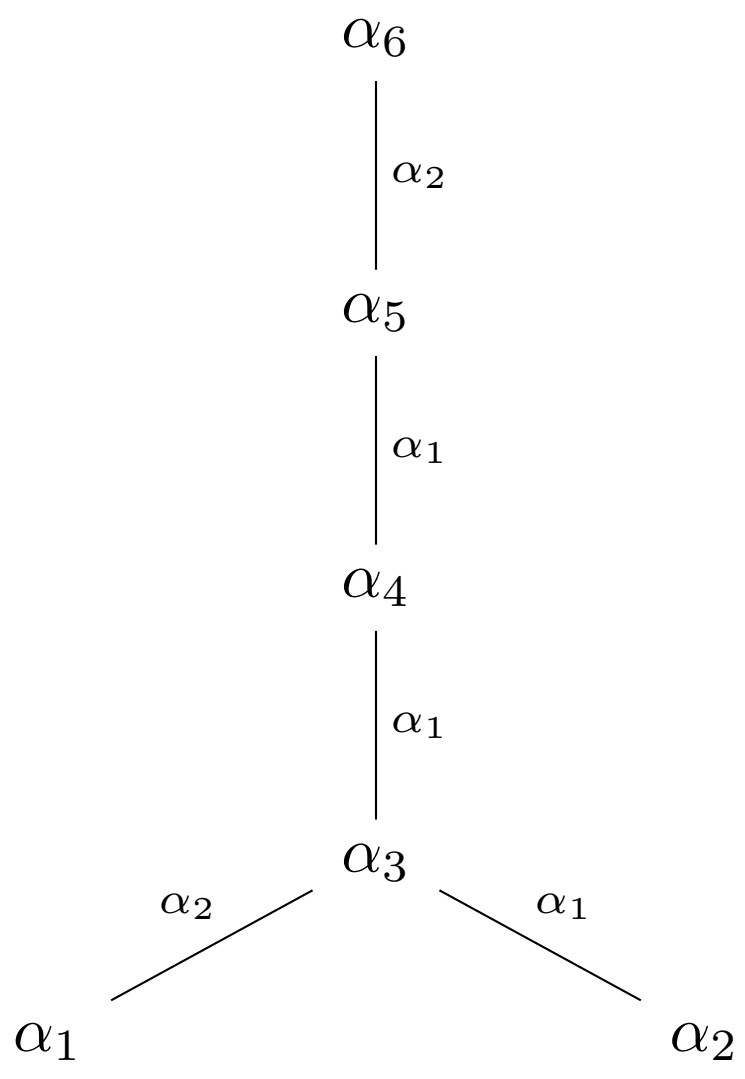
\begin{figure}[htbp]
  % https://tikzcd.yichuanshen.de/#N4Igdg9gJgpgziAXAbVABwnAlgFyxMJZAJgBoAGAXVJADcBDAGwFcYkQAdDptAC3oD6ANhABfUuky58hFGWLU6TVuy49+AgKxiJIDNjwEiZACyKGLNok7dGfQSZ2SDM46SHnlVm+sEBmJz0pQ1lkclIAdk9LVVt7AQBGQP1pIxQTSOiVazU7DWIxRRgoAHN4IlAAMwAnCABbJAyQHAgkPxpeGHoodkgwNhoLbJ88wQLxKtqGxE0aFraOrp7rPoGlGJy4jSSaRnoAIxhGAAVg12tqrBLeHECa+oXm1sQyEE7u3oI1oe9c+J2QHtDiczmkQJdrrcJiB7tNXvNEAlBl5Yr5EiBdgcjqcXGCITcMW8lp9+ncpkgkU8kOEiR8Vl9CT9UaMBAVMcCcalZOCrgTRJRREA
  \begin{tikzcd}[matrix xscale=0.8,matrix yscale=1,nodes={scale=1.2}]
                                            &  & \alpha_6                                  &  &                                            \\
                                            &  &                                           &  &                                            \\
                                            &  & \alpha_5 \arrow[uu, "\alpha_2"', no head] &  &                                            \\
                                            &  &                                           &  &                                            \\
                                            &  & \alpha_4 \arrow[uu, "\alpha_1"', no head] &  &                                            \\
                                            &  &                                           &  &                                            \\
                                            &  & \alpha_3 \arrow[uu, "\alpha_1"', no head] &  &                                            \\
  \alpha_1 \arrow[rru, pos=0.6, "\alpha_2", no head] &  &                                           &  & \alpha_2 \arrow[llu, pos=0.6, "\alpha_1"', no head]
  \end{tikzcd}
  \caption{Positive roots of $G_2$}
  \label{fig:G3-positive-roots}
\end{figure}

Retain the notation and choices of \Cref{subsec:structure_of_F4_and_G3}.
In the present case, $\mathfrak{a}_2\cong G_2$ and $\Pi_2=\{\alpha_1,\alpha_2\}$, where $\alpha_1$ is a short root and $\alpha_2$ is a long root.
The six positive roots of $\mathfrak{a}_2$ are shown in \Cref{fig:G3-positive-roots}.

Choose any nonzero elements $e_{\alpha_1}\in(\mathfrak{a}_2)_{\alpha_1}$ and $e_{\alpha_2}\in(\mathfrak{a}_2)_{\alpha_2}$.
For $i\in\{3,4,5,6\}$, define $e_{\alpha_i}\in(\mathfrak{a}_2)_{\alpha_i}$ as follows:
\[e_{\alpha_3}=[e_{\alpha_2},e_{\alpha_1}],\qquad e_{\alpha_4}=\frac{1}{2}[e_{\alpha_1},e_{\alpha_3}],\qquad e_{\alpha_5}=\frac{1}{3}[e_{\alpha_1},e_{\alpha_4}],\qquad e_{\alpha_6}=[e_{\alpha_2},e_{\alpha_5}].\]
For $i=1,\ldots,6$, write $h_i:=h_{\alpha_i}$ and let $e_{-\alpha_i}$ be the negative root vector determined by the convention in \Cref{subsec:structure_of_F4_and_G3}.
The first assertion of \Cref{prop:common_basis_construction} gives a basis of $\mathfrak{g}_{\bar{0}}$:
\[\{h,e_{\alpha},e_{-\alpha}\}\cup\left(\{h_1,h_2\}\cup\bigcup_{i=1}^{6}\{e_{\alpha_i},e_{-\alpha_i}\}\right).\]

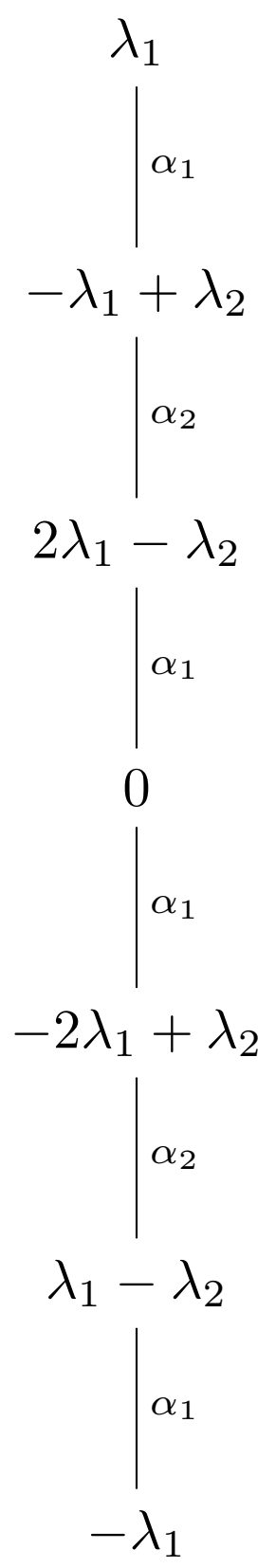
\begin{figure}[htbp]
  % https://tikzcd.yichuanshen.de/#N4Igdg9gJgpgziAXAbVABwnAlgFyxMJZABgBpiBdUkANwEMAbAVxiRAB12G6BbAIyh0A+gEYQAX1LpMufIRRkRVWoxZsAtJ279BogNRbeA4QCYJUkBmx4CRMieX1mrRCBOGdwkZq5HdZyWlrOTtSAGZHVRcQYnMg2VsFUgAWSOcNd19PfQ9jIQCLKwT5ElIAVjS1V1zdbxrTOMsZGxKyADZK6J9tPLFxZRgoAHN4IlAAMwAnCB4kMhAcCCQRagALGDooNkgwVmonKo52RjRVr0apmeXqRaQTNY2t1x29lXTq44ZThsCQS9nEPcFktEGF9lE2JwTmdRCAHpttgRWL9-kgwcCkMl4U9wEi4W9DlCvjC+hZUYgsRjEGVsYjdviDtEid98hdpgCaVS2rTnnjwe8jtDzv1xEA
  \begin{tikzcd}[matrix xscale=1,matrix yscale=1.8,nodes={scale=1.1}]
  \lambda_1 \arrow[d, "\alpha_1", no head]             \\
  -\lambda_1+\lambda_2 \arrow[d, "\alpha_2", no head]  \\
  2\lambda_1-\lambda_2 \arrow[d, "\alpha_1", no head]  \\
  0 \arrow[d, "\alpha_1", no head]                     \\
  -2\lambda_1+\lambda_2 \arrow[d, "\alpha_2", no head] \\
  \lambda_1-\lambda_2 \arrow[d, "\alpha_1", no head]   \\
  -\lambda_1                                          
  \end{tikzcd}
  \caption{Weights of $V_{G_2}(\lambda_1)$}
  \label{fig:G3-weight-diagram}
\end{figure}

We next specify the basis vectors in $\mathfrak{g}_{\bar{1}}\cong V_1\otimes V_2$.
The weights of $V_2$ are shown in \Cref{fig:G3-weight-diagram}.

\begin{rem}
  In terms of the simple roots, the fundamental weights are given by
  \[\lambda=\frac{1}{2}\alpha,\qquad\lambda_1=2\alpha_1+\alpha_2,\qquad\lambda_2=3\alpha_1+2\alpha_2.\]
\end{rem}

Choose any nonzero element $v_{-\lambda_1}\in (V_2)_{-\lambda_1}$.
In turn, for the other weights $\mu$ of $V_2$, define $v_{\mu}\in (V_2)_{\mu}$ as follows:
\[\begin{aligned}
  &v_{\lambda_1-\lambda_2}=e_{\alpha_1}\cdot v_{-\lambda_1}, &&\qquad v_{-2\lambda_1+\lambda_2}=e_{\alpha_2}\cdot v_{\lambda_1-\lambda_2}=e_{\alpha_3}\cdot v_{-\lambda_1},\\
  &v_0=e_{\alpha_1}\cdot v_{-2\lambda_1+\lambda_2}=e_{\alpha_4}\cdot v_{-\lambda_1}, &&\qquad v_{2\lambda_1-\lambda_2}=\frac{1}{2}e_{\alpha_1}\cdot v_0=e_{\alpha_5}\cdot v_{-\lambda_1},\\
  &v_{-\lambda_1+\lambda_2}=e_{\alpha_2}\cdot v_{2\lambda_1-\lambda_2}=e_{\alpha_6}\cdot v_{-\lambda_1}, &&\qquad v_{\lambda_1}=e_{\alpha_1}\cdot v_{-\lambda_1+\lambda_2}=e_{\alpha_6}\cdot v_{\lambda_1-\lambda_2}.
\end{aligned}\]

These vectors are nonzero by the root string facts recalled in \Cref{subsec:structure_of_F4_and_G3}.
The second assertion of \Cref{prop:common_basis_construction} therefore shows that, as $\mu$ ranges over the weights in \Cref{fig:G3-weight-diagram}, the $14$ elements $v_{\mu}^{+}$ and $v_{\mu}^{-}$ form a basis of $\mathfrak{g}_{\bar{1}}$.

\begin{rem}
  The roots $\lambda-\lambda_1$, $\alpha_1$, and $\alpha_2$ form a fundamental system of $G(3)$:
  \begin{equation}\label{eq:G3_fundamental_system}
    \Pi=\{\lambda-\lambda_1,\alpha_1,\alpha_2\}.
  \end{equation}
  Fix the order of the elements of $\Pi$ as shown in \eqref{eq:G3_fundamental_system}.
  The corresponding Cartan matrix, as recorded in \cite{Dictionary}, is
  \begin{equation}\label{eq:Cartan_mat_of_G(3)}
    \begin{pmatrix}0&1&0\\-1&2&-3\\0&-1&2\end{pmatrix}.
  \end{equation}
\end{rem}

Set $\delta=\lambda-\lambda_1$.
The vectors $v_{-\lambda_1}^{+}$ and $v_{\lambda_1}^{-}$ belong to the opposite odd root spaces $\mathfrak{g}_{\delta}$ and $\mathfrak{g}_{-\delta}$, respectively.
Standard facts about root spaces in basic simple Lie superalgebras imply that their bracket is a nonzero element of $\mathfrak{h}$; see \cite{CW}*{Theorem 1.18}.
The Cartan matrix \eqref{eq:Cartan_mat_of_G(3)} shows that this bracket is a nonzero scalar multiple of
\[2h+2h_1+3h_2.\]
Replacing $v_-$ by $cv_-$ also replaces $v_+$ by $cv_+$ and multiplies the bracket by $c^2$.
Since every nonzero complex number has a square root, we may therefore normalize $v_-$ so that
\begin{equation}\label{eq:G(3)_assumption}
  \left[v_{-\lambda_1}^{+},\,v_{\lambda_1}^{-}\right]=2h+2h_1+3h_2.
\end{equation}

The complete commutation relations among these basis elements are recorded in \Cref{sec:commutation_relations_for_G(3)}.

\subsection{Classification of Minimal \texorpdfstring{$\mathfrak{osp}(1\vert 2)$}{osp(1|2)}-Quintuples in \texorpdfstring{$G(3)$}{G(3)}}\label{subsec:G3_conjugacy}

By \Cref{subsec:minimal_nilpotent_orbits}, minimal $\mathfrak{sl}_2$-triples in $G(3)$ form two distinct conjugacy classes, represented by triples associated with $\alpha$ and with the highest root $\alpha_6$ of $\mathfrak{a}_2$.
Every minimal $\mathfrak{osp}(1\vert 2)$-quintuple extends a triple in one of these classes.
Using the diagonal inner automorphisms $\tau_{\beta,t}$ constructed in \Cref{subsec:inner_automorphisms} to rescale the root vectors, we may take these triples to be
\[\left(e_{\alpha},h,e_{-\alpha}\right)
  \qquad\text{or}\qquad
  \left(\frac{1}{3}e_{\alpha_6},h_6,3e_{-\alpha_6}\right).\]
Note that $h_6=h_1+2h_2$.
Conjugate quintuples have conjugate even triples, so we classify the extensions of each chosen triple under automorphisms fixing that triple.

For a quintuple extending the first triple, its fourth entry $y$ satisfies $[h,y]=-y$; hence $y$ lies in the span of the $v_{\mu}^{-}$ as $\mu$ ranges over the weights of $V_2$.

\begin{prop}\label{prop:G3_case_1_family}
  Let
  \[y=\sum_{\mu}b_{\mu}v_{\mu}^{-}\qquad(b_{\mu}\in\mathbb{C}),\]
  where the summation runs over the weights $\mu$ of $V_2$, and set
  \begin{equation}\label{eq:G3_case_1_scalar}
    b:=b_0^2-b_{-\lambda_1}b_{\lambda_1}
      +b_{\lambda_1-\lambda_2}b_{-\lambda_1+\lambda_2}
      -b_{-2\lambda_1+\lambda_2}b_{2\lambda_1-\lambda_2}.
  \end{equation}
  If $b\neq0$, set $X=e_{\alpha}/(4b)$, $Y=4be_{-\alpha}$, and $x=[X,y]$.
  Then $(X,x,h,y,Y)$ is a minimal $\mathfrak{osp}(1\vert 2)$-quintuple.
\end{prop}

\begin{proof}
  By the form of $y$, $[h,y]=-y$.  The relations in \Cref{sec:commutation_relations_for_G(3)} give
  \[-\frac{1}{2}[y,y]=4be_{-\alpha}.\]
  If $b\neq0$, the prescribed $(X,h,Y)$ is a minimal $\mathfrak{sl}_2$-triple by \Cref{prop:minimal_triple_entries} and \Cref{rem:scalar_multiple_minimal_triple}.
  These are the conditions in \Cref{lem:minimality_criterion}, so the assertion follows.
\end{proof}

In particular, the extensions of $(e_{\alpha},h,e_{-\alpha})$ are precisely the quintuples in \Cref{prop:G3_case_1_family} with $b=1/4$.

\begin{prop}\label{prop:G3_case_1_conjugacy}
  Suppose $(e_{\alpha},x_1,h,y_1,e_{-\alpha})$ and $(e_{\alpha},x_2,h,y_2,e_{-\alpha})$ are minimal $\mathfrak{osp}(1\vert 2)$-quintuples whose entries $y_1$ and $y_2$ are as in \Cref{prop:G3_case_1_family} with $b=1/4$.
  Then the two quintuples are conjugate.
\end{prop}

\begin{proof}
  Order the basis vectors in \Cref{prop:G3_case_1_family} as follows:
  \[\begin{gathered}
    w_0=v_{-\lambda_1}^{-},\quad
    w_1=v_{\lambda_1-\lambda_2}^{-},\quad
    w_2=v_{-2\lambda_1+\lambda_2}^{-},\quad
    w_3=v_0^{-},\\
    w_4=v_{2\lambda_1-\lambda_2}^{-},\quad
    w_5=v_{-\lambda_1+\lambda_2}^{-},\quad
    w_6=v_{\lambda_1}^{-}.
  \end{gathered}\]
  Write $y=\sum_{i=0}^6b_iw_i$.
  The condition $b=1/4$ becomes
  \begin{equation}\label{eq:G3_case_1_quadratic}
    b_3^2-b_0b_6+b_1b_5-b_2b_4=\frac{1}{4}.
  \end{equation}
  Set
  \[(\gamma_1,\gamma_2,\gamma_3,\gamma_4,\gamma_5)
    =(\alpha_1,\alpha_3,\alpha_4,\alpha_5,\alpha_6).\]
  The relations in \Cref{sec:commutation_relations_for_G(3)} give
  \[[e_{\gamma_i},w_0]=w_i
    \qquad\text{and}\qquad
    [e_{-\gamma_i},w_i]\in\mathbb{C}^{\times}w_0
    \qquad(1\leq i\leq5).\]
  Moreover, for $1\leq i\leq5$ and $1\leq j\leq6$, we have
  \[[e_{\gamma_i},w_j]\in\operatorname{span}_{\mathbb{C}}\{w_{i+1},\dotsc,w_6\}.\]

  Suppose first that $b_0=0$.
  Equation \eqref{eq:G3_case_1_quadratic} implies that $b_i\neq0$ for some $1\leq i\leq5$.
  For such an $i$, $[e_{-\gamma_i},w_i]\in\mathbb{C}^{\times}w_0$, and the distinct weights imply that no other $w_j$ contributes to the $w_0$-coefficient at first order.
  Thus the $w_0$-coefficient of $\sigma_{-\gamma_i,c}(y)$ is a polynomial in $c$ with a nonzero linear term.
  Choosing $c$ away from its zeros, we may assume $b_0\neq0$.
  Successively for $i=1,\dotsc,5$, apply $\sigma_{\gamma_i,c_i}$.
  At the $i$-th step, this automorphism changes the current coefficient of $w_i$ by $c_i b_0$ and leaves the coefficients of $w_0,\dotsc,w_{i-1}$ unchanged.
  We can therefore choose $c_i$ to cancel the coefficient of $w_i$ without reintroducing any coefficient canceled at an earlier step.
  We thereby obtain an element of the form
  \[\widetilde{y}=b_0w_0+\widetilde{b}_6w_6.\]
  All the automorphisms used fix $(e_{\alpha},h,e_{-\alpha})$ by \Cref{lem:root_exponentials_fix_triple}; hence \eqref{eq:G3_case_1_quadratic} remains valid and gives $\widetilde{b}_6=-1/(4b_0)$.
  Since $[h_1,w_0]=-w_0$ and $[h_1,w_6]=w_6$, the diagonal automorphism $\tau_{\alpha_1,2b_0}$ sends $\widetilde{y}$ to $(w_0-w_6)/2$.
  Thus both $y_1$ and $y_2$ are conjugate to $(w_0-w_6)/2$ by automorphisms fixing $(e_{\alpha},h,e_{-\alpha})$.
  It follows that the given two quintuples are conjugate.
\end{proof}

For a quintuple extending the second triple, its fourth entry $y$ satisfies $[h_6,y]=-y$; the weights in \Cref{fig:G3-weight-diagram} then force $y$ to lie in the span of $v_{-\lambda_1}^{\pm}$ and $v_{\lambda_1-\lambda_2}^{\pm}$.

\begin{prop}\label{prop:G3_case_2_family}
  Let $a_1,a_2,a_3,a_4\in\mathbb{C}$, and set
  \begin{equation}\label{eq:G3_case_2_element}
    y=a_1v_{-\lambda_1}^{-}+a_2v_{-\lambda_1}^{+}
      +a_3v_{\lambda_1-\lambda_2}^{-}+a_4v_{\lambda_1-\lambda_2}^{+}
    \qquad\text{and}\qquad
    a:=a_1a_4-a_2a_3.
  \end{equation}
  If $a\neq0$, set $X=e_{\alpha_6}/(3a)$, $Y=3ae_{-\alpha_6}$, and $x=[X,y]$.
  Then $(X,x,h_6,y,Y)$ is a minimal $\mathfrak{osp}(1\vert 2)$-quintuple.
\end{prop}

\begin{proof}
  By the form of $y$, $[h_6,y]=-y$.  The relations in \Cref{sec:commutation_relations_for_G(3)} give
  \[-\frac{1}{2}[y,y]=3ae_{-\alpha_6}.\]
  If $a\neq0$, the prescribed $(X,h_6,Y)$ is a minimal $\mathfrak{sl}_2$-triple by \Cref{prop:minimal_triple_entries} and \Cref{rem:scalar_multiple_minimal_triple}.
  These are the conditions in \Cref{lem:minimality_criterion}, so the assertion follows.
\end{proof}

In particular, the extensions of $(e_{\alpha_6}/3,h_6,3e_{-\alpha_6})$ are precisely the quintuples in \Cref{prop:G3_case_2_family} with $a=1$.

\begin{prop}\label{prop:G3_case_2_conjugacy}
  Suppose $(e_{\alpha_6}/3,x_1,h_6,y_1,3e_{-\alpha_6})$ and $(e_{\alpha_6}/3,x_2,h_6,y_2,3e_{-\alpha_6})$ are minimal $\mathfrak{osp}(1\vert 2)$-quintuples whose entries $y_1$ and $y_2$ are as in \eqref{eq:G3_case_2_element} with $a=1$.
  Then the two quintuples are conjugate.
\end{prop}

\begin{proof}
  Identify the element in \eqref{eq:G3_case_2_element} with the matrix
  \[A=\begin{pmatrix}a_1&a_2\\a_3&a_4\end{pmatrix}.\]
  The normalization $a=1$ means that $\det A=1$.
  Since $\alpha(h_6)=0$, \Cref{lem:root_exponentials_fix_triple} shows that $\sigma_{\alpha,c}$ and $\sigma_{-\alpha,c}$ fix the triple $(e_{\alpha_6}/3,h_6,3e_{-\alpha_6})$ pointwise.
  The relations in \Cref{sec:commutation_relations_for_G(3)} give
  \[\sigma_{\alpha,c}\cdot A
    =A\begin{pmatrix}1&c\\0&1\end{pmatrix},
    \qquad
    \sigma_{-\alpha,c}\cdot A
    =A\begin{pmatrix}1&0\\c&1\end{pmatrix}.\]
  These are right multiplications by the two types of elementary unipotent matrices, which generate $\operatorname{SL}_2(\mathbb{C})$.
  In particular, if $A_1,A_2\in\operatorname{SL}_2(\mathbb{C})$ correspond to $y_1,y_2$, then $A_1^{-1}A_2$ is a product of these elementary matrices, so $y_1$ and $y_2$ are conjugate by automorphisms fixing the triple.
  The same automorphisms therefore conjugate the quintuples.
\end{proof}

\subsection{The Centralizers and Complementary Submodules in \texorpdfstring{$G(3)$}{G(3)}}

For each conjugacy class found in the preceding subsection, we choose a representative $\mathcal{Q}$ with even triple $\mathcal{T}$ and compute the centralizers $\mathfrak{g}^{\mathcal{T}}$ and $\mathfrak{g}^{\mathcal{Q}}$.
We also construct a complementary submodule $M_{\mathcal{Q}}^{\mathfrak{g}}$ and determine its isomorphism type.

\begin{prop}\label{prop:G3_first_centralizers}
  Consider the minimal $\mathfrak{osp}(1\vert 2)$-quintuple $\mathcal{Q}$ given by
  \[\mathcal{Q}=\left(e_{\alpha},\frac{1}{2}v_0^{+},h,\frac{1}{2}v_0^{-},e_{-\alpha}\right),\]
  and write $\mathcal{T}$ for its even triple.
  Then
  \[\mathfrak{g}^{\mathcal{T}}=\mathfrak{a}_2\cong G_2,\qquad
    \mathfrak{g}^{\mathcal{Q}}\cong\mathfrak{sl}_3.\]
  Moreover, there exists a complementary submodule $M_{\mathcal{Q}}^{\mathfrak{g}}$, and
  \[M_{\mathcal{Q}}^{\mathfrak{g}}\cong
      V_{\mathfrak{sl}_3}(\varpi_1)\oplus V_{\mathfrak{sl}_3}(\varpi_2),\]
  where $\varpi_1$ and $\varpi_2$ are the fundamental weights of $\mathfrak{sl}_3$.
\end{prop}

\begin{proof}
  The entries of $\mathcal{T}$ span $\mathfrak{a}_1$, so their joint centralizer in $\mathfrak{a}_1$ is zero.
  Their joint centralizer in $\mathfrak{g}_{\bar{1}}\cong V_1\otimes V_2$ is also zero, because $V_1$ is the natural $\mathfrak{sl}_2$-module.
  Since $\mathfrak{a}_1$ commutes with $\mathfrak{a}_2$, we have $\mathfrak{g}^{\mathcal{T}}=\mathfrak{a}_2\cong G_2$.

  Write $y$ for the fourth entry of $\mathcal{Q}$.
  If
  \[g=\sum_{i=1}^{6}a_ie_{\alpha_i}
      +\sum_{i=1}^{6}b_ie_{-\alpha_i}
      +c_1h_1+c_2h_2\in\mathfrak{g}^{\mathcal{T}},\]
  then
  \[
    [g,y]=a_1v_{2\lambda_1-\lambda_2}^{-}
      +a_3v_{-\lambda_1+\lambda_2}^{-}
      +a_4v_{\lambda_1}^{-}
      +b_1v_{-2\lambda_1+\lambda_2}^{-}
      -b_3v_{\lambda_1-\lambda_2}^{-}
      +b_4v_{-\lambda_1}^{-}.
  \]
  It follows that
  \[\mathfrak{g}^{\mathcal{Q}}
    =\{g\in\mathfrak{g}^{\mathcal{T}}\mid[g,y]=0\}
    =\operatorname{span}_{\mathbb{C}}
      \{e_{\alpha_2},e_{\alpha_5},e_{\alpha_6},
        e_{-\alpha_2},e_{-\alpha_5},e_{-\alpha_6},h_1,h_2\}.\]
  The subalgebra $\mathfrak{h}_2=\operatorname{span}_{\mathbb{C}}\{h_1,h_2\}$ is a Cartan subalgebra of $\mathfrak{g}^{\mathcal{Q}}$.
  Relative to it, $\alpha_2$ and $\alpha_5$ form a fundamental system, with corresponding coroots $h_2$ and $h_1+h_2$.
  The relations in \Cref{sec:commutation_relations_for_G(3)} give the Cartan matrix
  \[\begin{pmatrix}2&-1\\-1&2\end{pmatrix}\]
  and $[e_{\alpha_2},e_{\alpha_5}]=e_{\alpha_6}$.
  Hence $\mathfrak{g}^{\mathcal{Q}}\cong\mathfrak{sl}_3$.
  Let $\varpi_1$ and $\varpi_2$ denote its fundamental weights in this order.

  Define
  \[M_{\mathcal{Q}}^{\mathfrak{g}}
    =\operatorname{span}_{\mathbb{C}}
      \{e_{\alpha_3},e_{\alpha_1},e_{-\alpha_4}\}
      \oplus\operatorname{span}_{\mathbb{C}}
      \{e_{\alpha_4},e_{-\alpha_1},e_{-\alpha_3}\}.\]
  The six listed vectors are linearly independent modulo $\mathfrak{g}^{\mathcal{Q}}$ and exhaust the complementary root spaces of $\mathfrak{g}^{\mathcal{T}}$.
  Hence $\mathfrak{g}^{\mathcal{T}}=\mathfrak{g}^{\mathcal{Q}}\oplus M_{\mathcal{Q}}^{\mathfrak{g}}$ as vector spaces.
  The relations in \Cref{sec:commutation_relations_for_G(3)} show that both summands of $M_{\mathcal{Q}}^{\mathfrak{g}}$ are preserved by $e_{\pm\alpha_2}$ and $e_{\pm\alpha_5}$.
  Since these four elements generate $\mathfrak{g}^{\mathcal{Q}}$, both summands are $\mathfrak{g}^{\mathcal{Q}}$-submodules.
  Relative to the ordered coroots $(h_2,h_1+h_2)$, the weights of the listed basis elements in the two summands are, respectively,
  \[\{(1,0),(-1,1),(0,-1)\}
    \qquad\text{and}\qquad
    \{(0,1),(1,-1),(-1,0)\}.\]
  Moreover, the corresponding simple root vectors $e_{\alpha_2}$ and $e_{\alpha_5}$ annihilate $e_{\alpha_3}$ and $e_{\alpha_4}$.
  Thus the two summands are isomorphic to $V_{\mathfrak{g}^{\mathcal{Q}}}(\varpi_1)$ and $V_{\mathfrak{g}^{\mathcal{Q}}}(\varpi_2)$, respectively.
\end{proof}

\begin{prop}\label{prop:G3_second_centralizers}
  Consider the minimal $\mathfrak{osp}(1\vert 2)$-quintuple $\mathcal{Q}$ given by
  \[\mathcal{Q}=\left(\frac{1}{3}e_{\alpha_6},
      \frac{1}{3}\bigl(v_{-\lambda_1+\lambda_2}^{-}+v_{\lambda_1}^{+}\bigr),
      h_6,
      v_{-\lambda_1}^{-}+v_{\lambda_1-\lambda_2}^{+},
      3e_{-\alpha_6}\right),\]
  and write $\mathcal{T}$ for its even triple.
  Then
  \[\mathfrak{g}^{\mathcal{T}}\cong\mathfrak{osp}(3\vert 2),\qquad
    \mathfrak{g}^{\mathcal{Q}}\cong\mathfrak{osp}(1\vert 2).\]
  Moreover, there exists a complementary submodule $M_{\mathcal{Q}}^{\mathfrak{g}}$, and
  \[M_{\mathcal{Q}}^{\mathfrak{g}}\cong\mathbb{C}^{3\vert 4}.\]
  Here, $\mathbb{C}^{3\vert 4}$ denotes the irreducible $\mathfrak{osp}(1\vert 2)$-module of dimension $(3\vert 4)$.
\end{prop}

\begin{proof}
  The bracket relations with $e_{\pm\alpha_6}$ give
  \[\mathfrak{g}^{\mathcal{T}}
    =\operatorname{span}_{\mathbb{C}}
      \{e_{\alpha},h,e_{-\alpha};
        e_{\alpha_1},h_1,e_{-\alpha_1}\}
      \oplus\operatorname{span}_{\mathbb{C}}
      \{v_{-2\lambda_1+\lambda_2}^{\pm},v_0^{\pm},
        v_{2\lambda_1-\lambda_2}^{\pm}\}.\]
  Its even part is the direct sum of the two commuting $\mathfrak{sl}_2$-subalgebras with respective bases
  \[\{e_{\alpha},h,e_{-\alpha}\}
    \qquad\text{and}\qquad
    \{e_{\alpha_1},h_1,e_{-\alpha_1}\}.\]
  Under these two summands, the odd part is the tensor product of the irreducible modules of dimensions $2$ and $3$, respectively, and hence is irreducible.
  Both simple summands act nontrivially on this tensor product, so it is also faithful as a module for the even part.
  The relations
  \[\left[v_{-2\lambda_1+\lambda_2}^{+},v_{2\lambda_1-\lambda_2}^{+}\right]
      =-4e_{\alpha},\qquad
    \left[v_{-2\lambda_1+\lambda_2}^{+},v_0^{-}\right]
      =2e_{-\alpha_1}\]
  give a nonzero root vector in each even simple ideal.
  By the super Jacobi identity, the span of odd-odd brackets is an ideal of $(\mathfrak{g}^{\mathcal{T}})_{\bar{0}}$.
  Since it has a nonzero component in each simple even ideal,
  \[\bigl[(\mathfrak{g}^{\mathcal{T}})_{\bar{1}},(\mathfrak{g}^{\mathcal{T}})_{\bar{1}}\bigr]=(\mathfrak{g}^{\mathcal{T}})_{\bar{0}}.\]
  Hence \Cref{lem:simplicity_from_odd_part} shows that $\mathfrak{g}^{\mathcal{T}}$ is simple.
  The classification of finite-dimensional simple Lie superalgebras therefore identifies it with $\mathfrak{osp}(3\vert 2)$.

  Write $y$ for the fourth entry of $\mathcal{Q}$.
  If
  \[g_{\bar{0}}=c_1e_{\alpha}+c_2e_{-\alpha}+c_3h
      +c_4e_{\alpha_1}+c_5e_{-\alpha_1}+c_6h_1
      \in(\mathfrak{g}^{\mathcal{T}})_{\bar{0}},\]
  then
  \[[g_{\bar{0}},y]
    =(c_3+c_6)v_{\lambda_1-\lambda_2}^{+}
      +(c_1+c_5)v_{-\lambda_1}^{+}
      +(c_2+c_4)v_{\lambda_1-\lambda_2}^{-}
      -(c_3+c_6)v_{-\lambda_1}^{-}.\]
  Thus,
  \[(\mathfrak{g}^{\mathcal{Q}})_{\bar{0}}
    =\left\{g_{\bar{0}}\in(\mathfrak{g}^{\mathcal{T}})_{\bar{0}}\mid[g_{\bar{0}},y]=0\right\}
    =\operatorname{span}_{\mathbb{C}}\{e_{\alpha}-e_{-\alpha_1},h-h_1,e_{-\alpha}-e_{\alpha_1}\}.\]
  Also, if
  \[g_{\bar{1}}=d_1v_{2\lambda_1-\lambda_2}^{+}
      +d_2v_{2\lambda_1-\lambda_2}^{-}
      +d_3v_0^{+}+d_4v_0^{-}
      +d_5v_{-2\lambda_1+\lambda_2}^{+}
      +d_6v_{-2\lambda_1+\lambda_2}^{-}
      \in(\mathfrak{g}^{\mathcal{T}})_{\bar{1}},\]
  then
  \[[g_{\bar{1}},y]
    =-3d_2e_{-\alpha_2}
      +(d_1-2d_4)e_{-\alpha_3}
      +(d_6-2d_3)e_{-\alpha_4}
      +3d_5e_{-\alpha_5}.\]
  The odd centralizer is determined by $d_2=d_5=0$, $d_1=2d_4$, and $d_6=2d_3$; thus
  \[(\mathfrak{g}^{\mathcal{Q}})_{\bar{1}}
    =\left\{g_{\bar{1}}\in(\mathfrak{g}^{\mathcal{T}})_{\bar{1}}\mid[g_{\bar{1}},y]=0\right\}
    =\operatorname{span}_{\mathbb{C}}\left\{\frac{1}{2}v_0^{+}+v_{-2\lambda_1+\lambda_2}^{-},
        v_{2\lambda_1-\lambda_2}^{+}+\frac{1}{2}v_0^{-}\right\}.\]
  Set
  \[\begin{gathered}
    X''=e_{\alpha}-e_{-\alpha_1},\qquad
    H''=h-h_1,\qquad
    Y''=e_{-\alpha}-e_{\alpha_1},\\
    x''=\frac{1}{2}v_0^{+}+v_{-2\lambda_1+\lambda_2}^{-},\qquad
    y''=-v_{2\lambda_1-\lambda_2}^{+}-\frac{1}{2}v_0^{-}.
  \end{gathered}\]
  The relations in \Cref{sec:commutation_relations_for_G(3)} show that $(X'',x'',H'',y'',Y'')$ is an $\mathfrak{osp}(1\vert 2)$-quintuple, and its entries span $\mathfrak{g}^{\mathcal{Q}}$ by the preceding calculations.
  Hence $\mathfrak{g}^{\mathcal{Q}}\cong\mathfrak{osp}(1\vert 2)$.

  Define
  \[\begin{aligned}
    M_{\mathcal{Q}}^{\mathfrak{g}}=\operatorname{span}_{\mathbb{C}}\{&
      v_{-2\lambda_1+\lambda_2}^{+},
      e_{-\alpha_1}-4e_{\alpha},
      v_0^{+}-v_{-2\lambda_1+\lambda_2}^{-},
      h_1-4h,\\
      &v_{2\lambda_1-\lambda_2}^{+}-v_0^{-},
      4e_{-\alpha}-e_{\alpha_1},
      -v_{2\lambda_1-\lambda_2}^{-}\}.
  \end{aligned}\]
  The seven vectors spanning $M_{\mathcal{Q}}^{\mathfrak{g}}$ lie in $\mathfrak{g}^{\mathcal{T}}$ and are linearly independent modulo $\mathfrak{g}^{\mathcal{Q}}$.
  Since $\dim\mathfrak{g}^{\mathcal{T}}=12$, $\dim\mathfrak{g}^{\mathcal{Q}}=5$, and $\dim M_{\mathcal{Q}}^{\mathfrak{g}}=7$, we obtain $\mathfrak{g}^{\mathcal{T}}=\mathfrak{g}^{\mathcal{Q}}\oplus M_{\mathcal{Q}}^{\mathfrak{g}}$ as vector spaces.
  Denote the seven listed vectors, in order, by
  \[m_3,m_2,m_1,m_0,m_{-1},m_{-2},m_{-3}.\]
  Direct calculation gives
  \[[H'',m_j]=jm_j\qquad(-3\leq j\leq3),\]
  together with
  \[\begin{aligned}
    {}[x'',m_2]&=3m_3,& [x'',m_1]&=-m_2,& [x'',m_0]&=2m_1,\\
    [x'',m_{-1}]&=-m_0,& [x'',m_{-2}]&=m_{-1},& [x'',m_{-3}]&=-m_{-2},
  \end{aligned}\]
  and
  \[\begin{aligned}
    {}[y'',m_3]&=-m_2,& [y'',m_2]&=-m_1,& [y'',m_1]&=-m_0,\\
    [y'',m_0]&=-2m_{-1},& [y'',m_{-1}]&=-m_{-2},& [y'',m_{-2}]&=-3m_{-3}.
  \end{aligned}\]
  The remaining actions of $x''$ and $y''$ on these vectors are zero.
  The elements $x''$ and $y''$ generate $\mathfrak{g}^{\mathcal{Q}}$, so these actions make $M_{\mathcal{Q}}^{\mathfrak{g}}$ a submodule.
  Each weight space is one-dimensional, every $m_j$ reaches $m_3$ under repeated application of $x''$, and repeated application of $y''$ to $m_3$ recovers all seven vectors.
  Hence the submodule is irreducible, with dimension $(3\vert 4)$, and is isomorphic to $\mathbb{C}^{3\vert 4}$.
\end{proof}

\begin{rem}
  The calculations in \Cref{prop:G3_first_centralizers,prop:G3_second_centralizers} show that, in each case, $\mathfrak{g}^{\mathcal{T}}$ is a direct sum of pairwise nonisomorphic irreducible $\mathfrak{g}^{\mathcal{Q}}$-submodules.
  Consequently, the complementary submodule $M_{\mathcal{Q}}^{\mathfrak{g}}$ constructed in each case is in fact the unique $\mathfrak{g}^{\mathcal{Q}}$-submodule of $\mathfrak{g}^{\mathcal{T}}$ complementary to $\mathfrak{g}^{\mathcal{Q}}$.
\end{rem}

\subsection{Summary for \texorpdfstring{$G(3)$}{G(3)}}

We summarize the results of this section.

\begin{thm}\label{thm:main_G(3)}
  Let $\mathfrak{g}=G(3)$, and let $\mathfrak{a}_1$ and $\mathfrak{a}_2$ be the simple ideals of $\mathfrak{g}_{\bar{0}}$ isomorphic to $\mathfrak{sl}_2$ and $G_2$, respectively.
  There are exactly two conjugacy classes of minimal $\mathfrak{osp}(1\vert 2)$-quintuples in $\mathfrak{g}$.
  Representatives $\mathcal{Q}_i=(X_i,x_i,H_i,y_i,Y_i)$, $i=1,2$, of these classes are given by
  \begin{align*}
    &X_1=e_{\alpha},
    &&x_1=\frac{1}{2}v_0^{+},
    &&H_1=h,\\
    &y_1=\frac{1}{2}v_0^{-},
    &&Y_1=e_{-\alpha},\\[2pt]
    &X_2=\frac{1}{3}e_{\alpha_6},
    &&x_2=\frac{1}{3}\bigl(v_{-\lambda_1+\lambda_2}^{-}+v_{\lambda_1}^{+}\bigr),
    &&H_2=h_6,\\
    &y_2=v_{-\lambda_1}^{-}+v_{\lambda_1-\lambda_2}^{+},
    &&Y_2=3e_{-\alpha_6}.
  \end{align*}
  Let $\mathcal{T}_i=(X_i,H_i,Y_i)$ be the even triple of $\mathcal{Q}_i$ for $i=1,2$.
  If $\varpi_1$ and $\varpi_2$ are the fundamental weights of $\mathfrak{sl}_3$, then the associated centralizers and complementary submodules are
  \begin{align*}
    \mathfrak{g}^{\mathcal{T}_1}&=\mathfrak{a}_2\cong G_2,
    &\mathfrak{g}^{\mathcal{Q}_1}&\cong\mathfrak{sl}_3,
    &M_{\mathcal{Q}_1}^{\mathfrak{g}}&\cong
      V_{\mathfrak{sl}_3}(\varpi_1)\oplus V_{\mathfrak{sl}_3}(\varpi_2),\\
    \mathfrak{g}^{\mathcal{T}_2}&\cong\mathfrak{osp}(3\vert 2),
    &\mathfrak{g}^{\mathcal{Q}_2}&\cong\mathfrak{osp}(1\vert 2),
    &M_{\mathcal{Q}_2}^{\mathfrak{g}}&\cong\mathbb{C}^{3\vert 4}.
  \end{align*}
  Here, $\mathbb{C}^{3\vert 4}$ denotes the irreducible $\mathfrak{osp}(1\vert 2)$-module of dimension $(3\vert 4)$.
\end{thm}

\begin{proof}
  The representatives and their centralizers and complementary submodules are given by \Cref{prop:G3_first_centralizers,prop:G3_second_centralizers}.
  The classification of minimal $\mathfrak{sl}_2$-triples in \Cref{subsec:minimal_nilpotent_orbits}, together with the normalization at the start of \Cref{subsec:G3_conjugacy}, shows that every minimal $\mathfrak{osp}(1\vert 2)$-quintuple is conjugate to an extension of $\mathcal{T}_1$ or $\mathcal{T}_2$.
  The extensions of each triple are conjugate by \Cref{prop:G3_case_1_conjugacy,prop:G3_case_2_conjugacy}; hence these representatives give exactly the two conjugacy classes.
\end{proof}

\appendix

\section{Commutation Relations}\label{sec:commutation_relations}

This appendix records all nonzero commutation relations, up to skew-supersymmetry, for the bases of $F(4)$ and $G(3)$ constructed in \Cref{subsec:explicit_basis_for_F(4),subsec:explicit_basis_for_G(3)}.
Whenever paired signs occur in a relation, either all upper signs or all lower signs are to be taken consistently throughout that relation.
All brackets not listed below and not determined by skew-supersymmetry are zero.

The relations follow from weight considerations, the normalizations and recursive definitions in the corresponding basis construction, and repeated applications of the super Jacobi identity.
Retain the notation of \Cref{subsec:structure_of_F4_and_G3}.
By definition, $[e_{-\alpha},v_{\mu}^{-}]=0$ and $[e_{\alpha},v_{\mu}^{-}]=v_{\mu}^{+}$.
Since $(e_{\alpha},h,e_{-\alpha})$ is an $\mathfrak{sl}_2$-triple, we also have
\[[e_{-\alpha},v_{\mu}^{+}]=[e_{-\alpha},[e_{\alpha},v_{\mu}^{-}]]
  =[[e_{-\alpha},e_{\alpha}],v_{\mu}^{-}]+[e_{\alpha},[e_{-\alpha},v_{\mu}^{-}]]
  =[-h,v_{\mu}^{-}]=v_{\mu}^{-}.\]
Thus the relations involving the common $\mathfrak{sl}_2$-factor are
\begin{equation}\label{eq:common_sl2_relations}
  [e_{\alpha},e_{-\alpha}]=h,\qquad
  [h,e_{\pm\alpha}]=\pm2e_{\pm\alpha},\qquad
  [e_{\mp\alpha},v_{\mu}^{\pm}]=v_{\mu}^{\mp},\qquad
  [h,v_{\mu}^{\pm}]=\pm v_{\mu}^{\pm}.
\end{equation}

The following observation reduces the number of odd-odd brackets that must be computed directly in both cases.
Let $\mu,\nu\in\mathfrak{h}_2^*$ be weights of $V_2$ with $\mu\neq-\nu$.
Then
\[[v_{\mu}^{\pm},v_{\nu}^{\pm}]=0,
  \qquad [v_{\mu}^{-},v_{\nu}^{+}]=-[v_{\mu}^{+},v_{\nu}^{-}].\]
Indeed, $[v_{\mu}^{+},v_{\nu}^{+}]\in\mathfrak{g}_{\alpha+\mu+\nu}$.
Since $\mu+\nu\neq0$ and every even root is supported on one of the two direct summands $\mathfrak{a}_1$ and $\mathfrak{a}_2$, the weight $\alpha+\mu+\nu$ is not a root; hence this bracket is zero.
The same argument gives $[v_{\mu}^{-},v_{\nu}^{-}]=0$.
Applying $\operatorname{ad}(e_{-\alpha})$ to $[v_{\mu}^{+},v_{\nu}^{+}]=0$ and using \eqref{eq:common_sl2_relations} together with the super Jacobi identity gives the second relation.

\subsection{Commutation Relations for \texorpdfstring{$F(4)$}{F(4)}}\label{sec:commutation_relations_for_F(4)}

We use the basis constructed in \Cref{subsec:explicit_basis_for_F(4)}.

The relations involving the common $\mathfrak{sl}_2$-factor are
\[[e_{\alpha},e_{-\alpha}]=h,\qquad
  [h,e_{\pm\alpha}]=\pm2e_{\pm\alpha},\qquad
  [e_{\mp\alpha},v_{\mu}^{\pm}]=v_{\mu}^{\mp},\qquad
  [h,v_{\mu}^{\pm}]=\pm v_{\mu}^{\pm}.\]

\paragraph{Even-even brackets.}
\[\begin{aligned}
  &[e_{\pm\alpha_1},e_{\pm\alpha_2}]=\mp e_{\pm\alpha_4}, &&[e_{\pm\alpha_1},e_{\pm\alpha_5}]=\pm e_{\pm\alpha_6}, &&[e_{\pm\alpha_1},e_{\pm\alpha_7}]=\pm e_{\pm\alpha_8},\\
  &[e_{\pm\alpha_2},e_{\pm\alpha_3}]=\pm e_{\pm\alpha_5}, &&[e_{\pm\alpha_2},e_{\pm\alpha_8}]=\pm e_{\pm\alpha_9}, &&[e_{\pm\alpha_3},e_{\pm\alpha_4}]=\pm e_{\pm\alpha_6},\\
  &[e_{\pm\alpha_3},e_{\pm\alpha_5}]=\pm2e_{\pm\alpha_7}, &&[e_{\pm\alpha_3},e_{\pm\alpha_6}]=\pm2e_{\pm\alpha_8}, &&[e_{\pm\alpha_4},e_{\pm\alpha_7}]=\pm e_{\pm\alpha_9},\\
  &[e_{\pm\alpha_5},e_{\pm\alpha_6}]=\pm2e_{\pm\alpha_9}.
\end{aligned}\]
\[\begin{aligned}
  &[e_{\alpha_1},e_{-\alpha_1}]=h_1, &&[e_{\pm\alpha_1},e_{\mp\alpha_4}]=\pm e_{\mp\alpha_2}, &&[e_{\pm\alpha_1},e_{\mp\alpha_6}]=\mp e_{\mp\alpha_5},\\
  &[e_{\pm\alpha_1},e_{\mp\alpha_8}]=\mp e_{\mp\alpha_7}, &&[e_{\alpha_2},e_{-\alpha_2}]=h_2, &&[e_{\pm\alpha_2},e_{\mp\alpha_4}]=\mp e_{\mp\alpha_1},\\
  &[e_{\pm\alpha_2},e_{\mp\alpha_5}]=\mp e_{\mp\alpha_3}, &&[e_{\pm\alpha_2},e_{\mp\alpha_9}]=\mp e_{\mp\alpha_8}, &&[e_{\alpha_3},e_{-\alpha_3}]=h_3,\\
  &[e_{\pm\alpha_3},e_{\mp\alpha_5}]=\pm2e_{\mp\alpha_2}, &&[e_{\pm\alpha_3},e_{\mp\alpha_6}]=\mp2e_{\mp\alpha_4}, &&[e_{\pm\alpha_3},e_{\mp\alpha_7}]=\mp e_{\mp\alpha_5},\\
  &[e_{\pm\alpha_3},e_{\mp\alpha_8}]=\mp e_{\mp\alpha_6}, &&[e_{\alpha_4},e_{-\alpha_4}]=h_1+h_2, &&[e_{\pm\alpha_4},e_{\mp\alpha_6}]=\pm e_{\mp\alpha_3},\\
  &[e_{\pm\alpha_4},e_{\mp\alpha_9}]=\mp e_{\mp\alpha_7}, &&[e_{\alpha_5},e_{-\alpha_5}]=2h_2+h_3, &&[e_{\pm\alpha_5},e_{\mp\alpha_6}]=\pm2e_{\mp\alpha_1},\\
  &[e_{\pm\alpha_5},e_{\mp\alpha_7}]=\pm e_{\mp\alpha_3}, &&[e_{\pm\alpha_5},e_{\mp\alpha_9}]=\mp e_{\mp\alpha_6}, &&[e_{\alpha_6},e_{-\alpha_6}]=2h_1+2h_2+h_3,\\
  &[e_{\pm\alpha_6},e_{\mp\alpha_8}]=\pm e_{\mp\alpha_3}, &&[e_{\pm\alpha_6},e_{\mp\alpha_9}]=\pm e_{\mp\alpha_5}, &&[e_{\alpha_7},e_{-\alpha_7}]=h_2+h_3,\\
  &[e_{\pm\alpha_7},e_{\mp\alpha_8}]=\pm e_{\mp\alpha_1}, &&[e_{\pm\alpha_7},e_{\mp\alpha_9}]=\pm e_{\mp\alpha_4}, &&[e_{\alpha_8},e_{-\alpha_8}]=h_1+h_2+h_3,\\
  &[e_{\pm\alpha_8},e_{\mp\alpha_9}]=\pm e_{\mp\alpha_2}, &&[e_{\alpha_9},e_{-\alpha_9}]=h_1+2h_2+h_3.
\end{aligned}\]
\[\begin{aligned}
  &[h_1,e_{\pm\alpha_1}]=\pm2e_{\pm\alpha_1}, &&[h_1,e_{\pm\alpha_2}]=\mp e_{\pm\alpha_2}, &&[h_1,e_{\pm\alpha_4}]=\pm e_{\pm\alpha_4},\\
  &[h_1,e_{\pm\alpha_5}]=\mp e_{\pm\alpha_5}, &&[h_1,e_{\pm\alpha_6}]=\pm e_{\pm\alpha_6}, &&[h_1,e_{\pm\alpha_7}]=\mp e_{\pm\alpha_7},\\
  &[h_1,e_{\pm\alpha_8}]=\pm e_{\pm\alpha_8}, &&[h_2,e_{\pm\alpha_1}]=\mp e_{\pm\alpha_1}, &&[h_2,e_{\pm\alpha_2}]=\pm2e_{\pm\alpha_2},\\
  &[h_2,e_{\pm\alpha_3}]=\mp e_{\pm\alpha_3}, &&[h_2,e_{\pm\alpha_4}]=\pm e_{\pm\alpha_4}, &&[h_2,e_{\pm\alpha_5}]=\pm e_{\pm\alpha_5},\\
  &[h_2,e_{\pm\alpha_8}]=\mp e_{\pm\alpha_8}, &&[h_2,e_{\pm\alpha_9}]=\pm e_{\pm\alpha_9}, &&[h_3,e_{\pm\alpha_2}]=\mp2e_{\pm\alpha_2},\\
  &[h_3,e_{\pm\alpha_3}]=\pm2e_{\pm\alpha_3}, &&[h_3,e_{\pm\alpha_4}]=\mp2e_{\pm\alpha_4}, &&[h_3,e_{\pm\alpha_7}]=\pm2e_{\pm\alpha_7},\\
  &[h_3,e_{\pm\alpha_8}]=\pm2e_{\pm\alpha_8}.
\end{aligned}\]

\paragraph{Even-odd brackets.}
\[\begin{aligned}
  &\left[e_{\alpha_1},v_{-\lambda_1+\lambda_2-\lambda_3}^{\pm}\right]=v_{\lambda_1-\lambda_3}^{\pm}, &&\quad\left[e_{\alpha_1},v_{-\lambda_1+\lambda_3}^{\pm}\right]=v_{\lambda_1-\lambda_2+\lambda_3}^{\pm}, &&\quad\left[e_{\alpha_2},v_{-\lambda_2+\lambda_3}^{\pm}\right]=v_{-\lambda_1+\lambda_2-\lambda_3}^{\pm},\\
  &\left[e_{\alpha_2},v_{\lambda_1-\lambda_2+\lambda_3}^{\pm}\right]=v_{\lambda_2-\lambda_3}^{\pm}, &&\quad\left[e_{\alpha_3},v_{-\lambda_3}^{\pm}\right]=v_{-\lambda_2+\lambda_3}^{\pm}, &&\quad\left[e_{\alpha_3},v_{-\lambda_1+\lambda_2-\lambda_3}^{\pm}\right]=v_{-\lambda_1+\lambda_3}^{\pm},\\
  &\left[e_{\alpha_3},v_{\lambda_1-\lambda_3}^{\pm}\right]=v_{\lambda_1-\lambda_2+\lambda_3}^{\pm}, &&\quad\left[e_{\alpha_3},v_{\lambda_2-\lambda_3}^{\pm}\right]=v_{\lambda_3}^{\pm}, &&\quad\left[e_{\alpha_4},v_{-\lambda_2+\lambda_3}^{\pm}\right]=-v_{\lambda_1-\lambda_3}^{\pm},\\
  &\left[e_{\alpha_4},v_{-\lambda_1+\lambda_3}^{\pm}\right]=v_{\lambda_2-\lambda_3}^{\pm}, &&\quad\left[e_{\alpha_5},v_{-\lambda_3}^{\pm}\right]=v_{-\lambda_1+\lambda_2-\lambda_3}^{\pm}, &&\quad\left[e_{\alpha_5},v_{-\lambda_2+\lambda_3}^{\pm}\right]=-v_{-\lambda_1+\lambda_3}^{\pm},\\
  &\left[e_{\alpha_5},v_{\lambda_1-\lambda_3}^{\pm}\right]=v_{\lambda_2-\lambda_3}^{\pm}, &&\quad\left[e_{\alpha_5},v_{\lambda_1-\lambda_2+\lambda_3}^{\pm}\right]=-v_{\lambda_3}^{\pm}, &&\quad\left[e_{\alpha_6},v_{-\lambda_3}^{\pm}\right]=v_{\lambda_1-\lambda_3}^{\pm},\\
  &\left[e_{\alpha_6},v_{-\lambda_2+\lambda_3}^{\pm}\right]=-v_{\lambda_1-\lambda_2+\lambda_3}^{\pm}, &&\quad\left[e_{\alpha_6},v_{-\lambda_1+\lambda_2-\lambda_3}^{\pm}\right]=-v_{\lambda_2-\lambda_3}^{\pm}, &&\quad\left[e_{\alpha_6},v_{-\lambda_1+\lambda_3}^{\pm}\right]=v_{\lambda_3}^{\pm},\\
  &\left[e_{\alpha_7},v_{-\lambda_3}^{\pm}\right]=v_{-\lambda_1+\lambda_3}^{\pm}, &&\quad\left[e_{\alpha_7},v_{\lambda_1-\lambda_3}^{\pm}\right]=v_{\lambda_3}^{\pm}, &&\quad\left[e_{\alpha_8},v_{-\lambda_3}^{\pm}\right]=v_{\lambda_1-\lambda_2+\lambda_3}^{\pm},\\
  &\left[e_{\alpha_8},v_{-\lambda_1+\lambda_2-\lambda_3}^{\pm}\right]=-v_{\lambda_3}^{\pm}, &&\quad\left[e_{\alpha_9},v_{-\lambda_3}^{\pm}\right]=v_{\lambda_2-\lambda_3}^{\pm}, &&\quad\left[e_{\alpha_9},v_{-\lambda_2+\lambda_3}^{\pm}\right]=v_{\lambda_3}^{\pm},
\end{aligned}\]
\[\begin{aligned}
  &\left[e_{-\alpha_1},v_{\lambda_1-\lambda_3}^{\pm}\right]=v_{-\lambda_1+\lambda_2-\lambda_3}^{\pm}, &&\left[e_{-\alpha_1},v_{\lambda_1-\lambda_2+\lambda_3}^{\pm}\right]=v_{-\lambda_1+\lambda_3}^{\pm}, &&\left[e_{-\alpha_2},v_{-\lambda_1+\lambda_2-\lambda_3}^{\pm}\right]=v_{-\lambda_2+\lambda_3}^{\pm},\\
  &\left[e_{-\alpha_2},v_{\lambda_2-\lambda_3}^{\pm}\right]=v_{\lambda_1-\lambda_2+\lambda_3}^{\pm}, &&\left[e_{-\alpha_3},v_{-\lambda_2+\lambda_3}^{\pm}\right]=v_{-\lambda_3}^{\pm}, &&\left[e_{-\alpha_3},v_{-\lambda_1+\lambda_3}^{\pm}\right]=v_{-\lambda_1+\lambda_2-\lambda_3}^{\pm},\\
  &\left[e_{-\alpha_3},v_{\lambda_1-\lambda_2+\lambda_3}^{\pm}\right]=v_{\lambda_1-\lambda_3}^{\pm}, &&\left[e_{-\alpha_3},v_{\lambda_3}^{\pm}\right]=v_{\lambda_2-\lambda_3}^{\pm}, &&\left[e_{-\alpha_4},v_{\lambda_1-\lambda_3}^{\pm}\right]=-v_{-\lambda_2+\lambda_3}^{\pm},\\
  &\left[e_{-\alpha_4},v_{\lambda_2-\lambda_3}^{\pm}\right]=v_{-\lambda_1+\lambda_3}^{\pm}, &&\left[e_{-\alpha_5},v_{-\lambda_1+\lambda_2-\lambda_3}^{\pm}\right]=v_{-\lambda_3}^{\pm}, &&\left[e_{-\alpha_5},v_{-\lambda_1+\lambda_3}^{\pm}\right]=-v_{-\lambda_2+\lambda_3}^{\pm},\\
  &\left[e_{-\alpha_5},v_{\lambda_2-\lambda_3}^{\pm}\right]=v_{\lambda_1-\lambda_3}^{\pm}, &&\left[e_{-\alpha_5},v_{\lambda_3}^{\pm}\right]=-v_{\lambda_1-\lambda_2+\lambda_3}^{\pm}, &&\left[e_{-\alpha_6},v_{\lambda_1-\lambda_3}^{\pm}\right]=v_{-\lambda_3}^{\pm},\\
  &\left[e_{-\alpha_6},v_{\lambda_1-\lambda_2+\lambda_3}^{\pm}\right]=-v_{-\lambda_2+\lambda_3}^{\pm}, &&\left[e_{-\alpha_6},v_{\lambda_2-\lambda_3}^{\pm}\right]=-v_{-\lambda_1+\lambda_2-\lambda_3}^{\pm}, &&\left[e_{-\alpha_6},v_{\lambda_3}^{\pm}\right]=v_{-\lambda_1+\lambda_3}^{\pm},\\
  &\left[e_{-\alpha_7},v_{-\lambda_1+\lambda_3}^{\pm}\right]=v_{-\lambda_3}^{\pm}, &&\left[e_{-\alpha_7},v_{\lambda_3}^{\pm}\right]=v_{\lambda_1-\lambda_3}^{\pm}, &&\left[e_{-\alpha_8},v_{\lambda_1-\lambda_2+\lambda_3}^{\pm}\right]=v_{-\lambda_3}^{\pm},\\
  &\left[e_{-\alpha_8},v_{\lambda_3}^{\pm}\right]=-v_{-\lambda_1+\lambda_2-\lambda_3}^{\pm}, &&\left[e_{-\alpha_9},v_{\lambda_2-\lambda_3}^{\pm}\right]=v_{-\lambda_3}^{\pm}, &&\left[e_{-\alpha_9},v_{\lambda_3}^{\pm}\right]=v_{-\lambda_2+\lambda_3}^{\pm},
\end{aligned}\]
\[\begin{aligned}
  &\left[h_1,v_{-\lambda_1+\lambda_2-\lambda_3}^{\pm}\right]=-v_{-\lambda_1+\lambda_2-\lambda_3}^{\pm}, &&\left[h_1,v_{\lambda_1-\lambda_3}^{\pm}\right]=v_{\lambda_1-\lambda_3}^{\pm}, &&\left[h_1,v_{-\lambda_1+\lambda_3}^{\pm}\right]=-v_{-\lambda_1+\lambda_3}^{\pm},\\
  &\left[h_1,v_{\lambda_1-\lambda_2+\lambda_3}^{\pm}\right]=v_{\lambda_1-\lambda_2+\lambda_3}^{\pm}, &&\left[h_2,v_{-\lambda_2+\lambda_3}^{\pm}\right]=-v_{-\lambda_2+\lambda_3}^{\pm}, &&\left[h_2,v_{-\lambda_1+\lambda_2-\lambda_3}^{\pm}\right]=v_{-\lambda_1+\lambda_2-\lambda_3}^{\pm},\\
  &\left[h_2,v_{\lambda_1-\lambda_2+\lambda_3}^{\pm}\right]=-v_{\lambda_1-\lambda_2+\lambda_3}^{\pm}, &&\left[h_2,v_{\lambda_2-\lambda_3}^{\pm}\right]=v_{\lambda_2-\lambda_3}^{\pm}, &&\left[h_3,v_{-\lambda_3}^{\pm}\right]=-v_{-\lambda_3}^{\pm},\\
  &\left[h_3,v_{-\lambda_2+\lambda_3}^{\pm}\right]=v_{-\lambda_2+\lambda_3}^{\pm}, &&\left[h_3,v_{-\lambda_1+\lambda_2-\lambda_3}^{\pm}\right]=-v_{-\lambda_1+\lambda_2-\lambda_3}^{\pm}, &&\left[h_3,v_{\lambda_1-\lambda_3}^{\pm}\right]=-v_{\lambda_1-\lambda_3}^{\pm},\\
  &\left[h_3,v_{-\lambda_1+\lambda_3}^{\pm}\right]=v_{-\lambda_1+\lambda_3}^{\pm}, &&\left[h_3,v_{\lambda_1-\lambda_2+\lambda_3}^{\pm}\right]=v_{\lambda_1-\lambda_2+\lambda_3}^{\pm}, &&\left[h_3,v_{\lambda_2-\lambda_3}^{\pm}\right]=-v_{\lambda_2-\lambda_3}^{\pm},\\
  &\left[h_3,v_{\lambda_3}^{\pm}\right]=v_{\lambda_3}^{\pm},
\end{aligned}\]

\paragraph{Odd-odd brackets.}
\[\begin{aligned}
  &\left[v_{-\lambda_3}^{\pm},\,v_{\lambda_3}^{\pm}\right]=\mp 3e_{\pm\alpha}, &&\qquad\left[v_{-\lambda_2+\lambda_3}^{\pm},\,v_{\lambda_2-\lambda_3}^{\pm}\right]=\pm 3e_{\pm\alpha},\\
  &\left[v_{-\lambda_1+\lambda_2-\lambda_3}^{\pm},\,v_{\lambda_1-\lambda_2+\lambda_3}^{\pm}\right]=\mp 3e_{\pm\alpha}, &&\qquad\left[v_{\lambda_1-\lambda_3}^{\pm},\,v_{-\lambda_1+\lambda_3}^{\pm}\right]=\pm 3e_{\pm\alpha},
\end{aligned}\]
\[\begin{aligned}
  &\left[v_{-\lambda_3}^{\pm},\,v_{-\lambda_2+\lambda_3}^{\mp}\right]=\pm 2e_{-\alpha_9}, &&\left[v_{-\lambda_3}^{\pm},\,v_{-\lambda_1+\lambda_2-\lambda_3}^{\mp}\right]=\mp 2e_{-\alpha_8},\\
  &\left[v_{-\lambda_3}^{\pm},\,v_{\lambda_1-\lambda_3}^{\mp}\right]=\pm 2e_{-\alpha_7}, &&\left[v_{-\lambda_3}^{\pm},\,v_{-\lambda_1+\lambda_3}^{\mp}\right]=\pm e_{-\alpha_6},\\
  &\left[v_{-\lambda_3}^{\pm},\,v_{\lambda_1-\lambda_2+\lambda_3}^{\mp}\right]=\mp e_{-\alpha_5}, &&\left[v_{-\lambda_3}^{\pm},\,v_{\lambda_2-\lambda_3}^{\mp}\right]=\pm e_{-\alpha_3},\\
  &\left[v_{-\lambda_3}^{\pm},\,v_{\lambda_3}^{\mp}\right]=\frac{3}{2}h\pm\left(h_1+2h_2+\frac{3}{2}h_3\right), &&\left[v_{-\lambda_2+\lambda_3}^{\pm},\,v_{-\lambda_1+\lambda_2-\lambda_3}^{\mp}\right]=\pm e_{-\alpha_6},\\
  &\left[v_{-\lambda_2+\lambda_3}^{\pm},\,v_{\lambda_1-\lambda_3}^{\mp}\right]=\mp e_{-\alpha_5}, &&\left[v_{-\lambda_2+\lambda_3}^{\pm},\,v_{-\lambda_1+\lambda_3}^{\mp}\right]=\mp 2e_{-\alpha_4},\\
  &\left[v_{-\lambda_2+\lambda_3}^{\pm},\,v_{\lambda_1-\lambda_2+\lambda_3}^{\mp}\right]=\mp 2e_{-\alpha_2}, &&\left[v_{-\lambda_2+\lambda_3}^{\pm},\,v_{\lambda_2-\lambda_3}^{\mp}\right]=-\frac{3}{2}h\mp\left(h_1+2h_2+\frac{1}{2}h_3\right),\\
  &\left[v_{-\lambda_2+\lambda_3}^{\pm},\,v_{\lambda_3}^{\mp}\right]=\mp e_{\alpha_3}, &&\left[v_{-\lambda_1+\lambda_2-\lambda_3}^{\pm},\,v_{\lambda_1-\lambda_3}^{\mp}\right]=\pm e_{-\alpha_3},\\
  &\left[v_{-\lambda_1+\lambda_2-\lambda_3}^{\pm},\,v_{-\lambda_1+\lambda_3}^{\mp}\right]=\pm 2e_{-\alpha_1}, &&\left[v_{-\lambda_1+\lambda_2-\lambda_3}^{\pm},\,v_{\lambda_1-\lambda_2+\lambda_3}^{\mp}\right]=\frac{3}{2}h\pm\left(h_1+\frac{1}{2}h_3\right),\\
  &\left[v_{-\lambda_1+\lambda_2-\lambda_3}^{\pm},\,v_{\lambda_2-\lambda_3}^{\mp}\right]=\pm 2e_{\alpha_2}, &&\left[v_{-\lambda_1+\lambda_2-\lambda_3}^{\pm},\,v_{\lambda_3}^{\mp}\right]=\mp e_{\alpha_5},\\
  &\left[v_{\lambda_1-\lambda_3}^{\pm},\,v_{-\lambda_1+\lambda_3}^{\mp}\right]=-\frac{3}{2}h\pm\left(h_1-\frac{1}{2}h_3\right), &&\left[v_{\lambda_1-\lambda_3}^{\pm},\,v_{\lambda_1-\lambda_2+\lambda_3}^{\mp}\right]=\mp 2e_{\alpha_1},\\
  &\left[v_{\lambda_1-\lambda_3}^{\pm},\,v_{\lambda_2-\lambda_3}^{\mp}\right]=\mp 2e_{\alpha_4}, &&\left[v_{\lambda_1-\lambda_3}^{\pm},\,v_{\lambda_3}^{\mp}\right]=\mp e_{\alpha_6},\\
  &\left[v_{-\lambda_1+\lambda_3}^{\pm},\,v_{\lambda_1-\lambda_2+\lambda_3}^{\mp}\right]=\mp e_{\alpha_3}, &&\left[v_{-\lambda_1+\lambda_3}^{\pm},\,v_{\lambda_2-\lambda_3}^{\mp}\right]=\mp e_{\alpha_5},\\
  &\left[v_{-\lambda_1+\lambda_3}^{\pm},\,v_{\lambda_3}^{\mp}\right]=\mp 2e_{\alpha_7}, &&\left[v_{\lambda_1-\lambda_2+\lambda_3}^{\pm},\,v_{\lambda_2-\lambda_3}^{\mp}\right]=\mp e_{\alpha_6},\\
  &\left[v_{\lambda_1-\lambda_2+\lambda_3}^{\pm},\,v_{\lambda_3}^{\mp}\right]=\mp 2e_{\alpha_8}, &&\left[v_{\lambda_2-\lambda_3}^{\pm},\,v_{\lambda_3}^{\mp}\right]=\mp 2e_{\alpha_9}.
\end{aligned}\]

\subsection{Commutation Relations for \texorpdfstring{$G(3)$}{G(3)}}\label{sec:commutation_relations_for_G(3)}

We use the basis constructed in \Cref{subsec:explicit_basis_for_G(3)}.

The relations involving the common $\mathfrak{sl}_2$-factor are
\[[e_{\alpha},e_{-\alpha}]=h,\qquad
  [h,e_{\pm\alpha}]=\pm2e_{\pm\alpha},\qquad
  [e_{\mp\alpha},v_{\mu}^{\pm}]=v_{\mu}^{\mp},\qquad
  [h,v_{\mu}^{\pm}]=\pm v_{\mu}^{\pm}.\]

\paragraph{Even-even brackets.}
\[\begin{aligned}
  &[e_{\pm\alpha_1},e_{\pm\alpha_2}]=\mp e_{\pm\alpha_3}, &&[e_{\pm\alpha_1},e_{\pm\alpha_3}]=\pm2e_{\pm\alpha_4}, &&[e_{\pm\alpha_1},e_{\pm\alpha_4}]=\pm3e_{\pm\alpha_5},\\
  &[e_{\pm\alpha_2},e_{\pm\alpha_5}]=\pm e_{\pm\alpha_6}, &&[e_{\pm\alpha_3},e_{\pm\alpha_4}]=\pm3e_{\pm\alpha_6}.
\end{aligned}\]
\[\begin{aligned}
  &[e_{\alpha_1},e_{-\alpha_1}]=h_1, &&[e_{\pm\alpha_1},e_{\mp\alpha_3}]=\pm3e_{\mp\alpha_2}, &&[e_{\pm\alpha_1},e_{\mp\alpha_4}]=\mp2e_{\mp\alpha_3},\\
  &[e_{\pm\alpha_1},e_{\mp\alpha_5}]=\mp e_{\mp\alpha_4}, &&[e_{\alpha_2},e_{-\alpha_2}]=h_2, &&[e_{\pm\alpha_2},e_{\mp\alpha_3}]=\mp e_{\mp\alpha_1},\\
  &[e_{\pm\alpha_2},e_{\mp\alpha_6}]=\mp e_{\mp\alpha_5}, &&[e_{\alpha_3},e_{-\alpha_3}]=h_1+3h_2, &&[e_{\pm\alpha_3},e_{\mp\alpha_4}]=\pm2e_{\mp\alpha_1},\\
  &[e_{\pm\alpha_3},e_{\mp\alpha_6}]=\mp e_{\mp\alpha_4}, &&[e_{\alpha_4},e_{-\alpha_4}]=2h_1+3h_2, &&[e_{\pm\alpha_4},e_{\mp\alpha_5}]=\pm e_{\mp\alpha_1},\\
  &[e_{\pm\alpha_4},e_{\mp\alpha_6}]=\pm e_{\mp\alpha_3}, &&[e_{\alpha_5},e_{-\alpha_5}]=h_1+h_2, &&[e_{\pm\alpha_5},e_{\mp\alpha_6}]=\pm e_{\mp\alpha_2},\\
  &[e_{\alpha_6},e_{-\alpha_6}]=h_1+2h_2.
\end{aligned}\]
\[\begin{aligned}
  &[h_1,e_{\pm\alpha_1}]=\pm2e_{\pm\alpha_1}, &&[h_1,e_{\pm\alpha_2}]=\mp3e_{\pm\alpha_2}, &&[h_1,e_{\pm\alpha_3}]=\mp e_{\pm\alpha_3},\\
  &[h_1,e_{\pm\alpha_4}]=\pm e_{\pm\alpha_4}, &&[h_1,e_{\pm\alpha_5}]=\pm3e_{\pm\alpha_5}, &&[h_2,e_{\pm\alpha_1}]=\mp e_{\pm\alpha_1},\\
  &[h_2,e_{\pm\alpha_2}]=\pm2e_{\pm\alpha_2}, &&[h_2,e_{\pm\alpha_3}]=\pm e_{\pm\alpha_3}, &&[h_2,e_{\pm\alpha_5}]=\mp e_{\pm\alpha_5},\\
  &[h_2,e_{\pm\alpha_6}]=\pm e_{\pm\alpha_6}.
\end{aligned}\]

\paragraph{Even-odd brackets.}
\[\begin{aligned}
  &\left[e_{\alpha_1},v_{-\lambda_1}^{\pm}\right]=v_{\lambda_1-\lambda_2}^{\pm}, &&\qquad\left[e_{\alpha_1},v_{-2\lambda_1+\lambda_2}^{\pm}\right]=v_{0}^{\pm}, &&\qquad\left[e_{\alpha_1},v_{0}^{\pm}\right]=2v_{2\lambda_1-\lambda_2}^{\pm},\\
  &\left[e_{\alpha_1},v_{-\lambda_1+\lambda_2}^{\pm}\right]=v_{\lambda_1}^{\pm}, &&\qquad\left[e_{\alpha_2},v_{\lambda_1-\lambda_2}^{\pm}\right]=v_{-2\lambda_1+\lambda_2}^{\pm}, &&\qquad\left[e_{\alpha_2},v_{2\lambda_1-\lambda_2}^{\pm}\right]=v_{-\lambda_1+\lambda_2}^{\pm},\\
  &\left[e_{\alpha_3},v_{-\lambda_1}^{\pm}\right]=v_{-2\lambda_1+\lambda_2}^{\pm}, &&\qquad\left[e_{\alpha_3},v_{\lambda_1-\lambda_2}^{\pm}\right]=-v_{0}^{\pm}, &&\qquad\left[e_{\alpha_3},v_{0}^{\pm}\right]=2v_{-\lambda_1+\lambda_2}^{\pm},\\
  &\left[e_{\alpha_3},v_{2\lambda_1-\lambda_2}^{\pm}\right]=-v_{\lambda_1}^{\pm}, &&\qquad\left[e_{\alpha_4},v_{-\lambda_1}^{\pm}\right]=v_{0}^{\pm}, &&\qquad\left[e_{\alpha_4},v_{\lambda_1-\lambda_2}^{\pm}\right]=-v_{2\lambda_1-\lambda_2}^{\pm},\\
  &\left[e_{\alpha_4},v_{-2\lambda_1+\lambda_2}^{\pm}\right]=-v_{-\lambda_1+\lambda_2}^{\pm}, &&\qquad\left[e_{\alpha_4},v_{0}^{\pm}\right]=2v_{\lambda_1}^{\pm}, &&\qquad\left[e_{\alpha_5},v_{-\lambda_1}^{\pm}\right]=v_{2\lambda_1-\lambda_2}^{\pm},\\
  &\left[e_{\alpha_5},v_{-2\lambda_1+\lambda_2}^{\pm}\right]=-v_{\lambda_1}^{\pm}, &&\qquad\left[e_{\alpha_6},v_{-\lambda_1}^{\pm}\right]=v_{-\lambda_1+\lambda_2}^{\pm}, &&\qquad\left[e_{\alpha_6},v_{\lambda_1-\lambda_2}^{\pm}\right]=v_{\lambda_1}^{\pm},
\end{aligned}\]
\[\begin{aligned}
  &\left[e_{-\alpha_1},v_{\lambda_1-\lambda_2}^{\pm}\right]=v_{-\lambda_1}^{\pm}, &&\qquad\left[e_{-\alpha_1},v_{0}^{\pm}\right]=2v_{-2\lambda_1+\lambda_2}^{\pm}, &&\qquad\left[e_{-\alpha_1},v_{2\lambda_1-\lambda_2}^{\pm}\right]=v_{0}^{\pm},\\
  &\left[e_{-\alpha_1},v_{\lambda_1}^{\pm}\right]=v_{-\lambda_1+\lambda_2}^{\pm}, &&\qquad\left[e_{-\alpha_2},v_{-2\lambda_1+\lambda_2}^{\pm}\right]=v_{\lambda_1-\lambda_2}^{\pm}, &&\qquad\left[e_{-\alpha_2},v_{-\lambda_1+\lambda_2}^{\pm}\right]=v_{2\lambda_1-\lambda_2}^{\pm},\\
  &\left[e_{-\alpha_3},v_{-2\lambda_1+\lambda_2}^{\pm}\right]=v_{-\lambda_1}^{\pm}, &&\qquad\left[e_{-\alpha_3},v_{0}^{\pm}\right]=-2v_{\lambda_1-\lambda_2}^{\pm}, &&\qquad\left[e_{-\alpha_3},v_{-\lambda_1+\lambda_2}^{\pm}\right]=v_{0}^{\pm},\\
  &\left[e_{-\alpha_3},v_{\lambda_1}^{\pm}\right]=-v_{2\lambda_1-\lambda_2}^{\pm}, &&\qquad\left[e_{-\alpha_4},v_{0}^{\pm}\right]=2v_{-\lambda_1}^{\pm}, &&\qquad\left[e_{-\alpha_4},v_{2\lambda_1-\lambda_2}^{\pm}\right]=-v_{\lambda_1-\lambda_2}^{\pm},\\
  &\left[e_{-\alpha_4},v_{-\lambda_1+\lambda_2}^{\pm}\right]=-v_{-2\lambda_1+\lambda_2}^{\pm}, &&\qquad\left[e_{-\alpha_4},v_{\lambda_1}^{\pm}\right]=v_{0}^{\pm}, &&\qquad\left[e_{-\alpha_5},v_{2\lambda_1-\lambda_2}^{\pm}\right]=v_{-\lambda_1}^{\pm},\\
  &\left[e_{-\alpha_5},v_{\lambda_1}^{\pm}\right]=-v_{-2\lambda_1+\lambda_2}^{\pm}, &&\qquad\left[e_{-\alpha_6},v_{-\lambda_1+\lambda_2}^{\pm}\right]=v_{-\lambda_1}^{\pm}, &&\qquad\left[e_{-\alpha_6},v_{\lambda_1}^{\pm}\right]=v_{\lambda_1-\lambda_2}^{\pm},
\end{aligned}\]
\[\begin{aligned}
  &\left[h_1,v_{-\lambda_1}^{\pm}\right]=-v_{-\lambda_1}^{\pm}, &&\qquad\left[h_1,v_{\lambda_1-\lambda_2}^{\pm}\right]=v_{\lambda_1-\lambda_2}^{\pm}, &&\qquad\left[h_1,v_{-2\lambda_1+\lambda_2}^{\pm}\right]=-2v_{-2\lambda_1+\lambda_2}^{\pm},\\
  &\left[h_1,v_{2\lambda_1-\lambda_2}^{\pm}\right]=2v_{2\lambda_1-\lambda_2}^{\pm}, &&\qquad\left[h_1,v_{-\lambda_1+\lambda_2}^{\pm}\right]=-v_{-\lambda_1+\lambda_2}^{\pm}, &&\qquad\left[h_1,v_{\lambda_1}^{\pm}\right]=v_{\lambda_1}^{\pm},\\
  &\left[h_2,v_{\lambda_1-\lambda_2}^{\pm}\right]=-v_{\lambda_1-\lambda_2}^{\pm}, &&\qquad\left[h_2,v_{-2\lambda_1+\lambda_2}^{\pm}\right]=v_{-2\lambda_1+\lambda_2}^{\pm}, &&\qquad\left[h_2,v_{2\lambda_1-\lambda_2}^{\pm}\right]=-v_{2\lambda_1-\lambda_2}^{\pm},\\
  &\left[h_2,v_{-\lambda_1+\lambda_2}^{\pm}\right]=v_{-\lambda_1+\lambda_2}^{\pm},
\end{aligned}\]

\paragraph{Odd-odd brackets.}
\[\begin{aligned}
  &\left[v_{-\lambda_1}^{\pm},\,v_{\lambda_1}^{\pm}\right]=\mp 4e_{\pm\alpha}, &&\qquad\left[v_{\lambda_1-\lambda_2}^{\pm},\,v_{-\lambda_1+\lambda_2}^{\pm}\right]=\pm 4e_{\pm\alpha},\\
  &\left[v_{-2\lambda_1+\lambda_2}^{\pm},\,v_{2\lambda_1-\lambda_2}^{\pm}\right]=\mp 4e_{\pm\alpha}, &&\qquad\left[v_{0}^{\pm},\,v_{0}^{\pm}\right]=\pm 8e_{\pm\alpha},
\end{aligned}\]
\[\begin{aligned}
  &\left[v_{-\lambda_1}^{\pm},\,v_{\lambda_1-\lambda_2}^{\mp}\right]=\pm 3e_{-\alpha_6}, &&\left[v_{-\lambda_1}^{\pm},\,v_{-2\lambda_1+\lambda_2}^{\mp}\right]=\mp 3e_{-\alpha_5}, &&\left[v_{-\lambda_1}^{\pm},\,v_{0}^{\mp}\right]=\pm 2e_{-\alpha_4},\\
  &\left[v_{-\lambda_1}^{\pm},\,v_{2\lambda_1-\lambda_2}^{\mp}\right]=\mp e_{-\alpha_3}, &&\left[v_{-\lambda_1}^{\pm},\,v_{-\lambda_1+\lambda_2}^{\mp}\right]=\pm e_{-\alpha_1}, &&\left[v_{-\lambda_1}^{\pm},\,v_{\lambda_1}^{\mp}\right]=2h\pm(2h_1+3h_2),\\
  &\left[v_{\lambda_1-\lambda_2}^{\pm},\,v_{-2\lambda_1+\lambda_2}^{\mp}\right]=\pm e_{-\alpha_4}, &&\left[v_{\lambda_1-\lambda_2}^{\pm},\,v_{0}^{\mp}\right]=\mp 2e_{-\alpha_3}, &&\left[v_{\lambda_1-\lambda_2}^{\pm},\,v_{2\lambda_1-\lambda_2}^{\mp}\right]=\mp 3e_{-\alpha_2},\\
  &\left[v_{\lambda_1-\lambda_2}^{\pm},\,v_{-\lambda_1+\lambda_2}^{\mp}\right]=-2h\mp(h_1+3h_2), &&\left[v_{\lambda_1-\lambda_2}^{\pm},\,v_{\lambda_1}^{\mp}\right]=\mp e_{\alpha_1}, &&\left[v_{-2\lambda_1+\lambda_2}^{\pm},\,v_{0}^{\mp}\right]=\pm 2e_{-\alpha_1},\\
  &\left[v_{-2\lambda_1+\lambda_2}^{\pm},\,v_{2\lambda_1-\lambda_2}^{\mp}\right]=2h\pm h_1, &&\left[v_{-2\lambda_1+\lambda_2}^{\pm},\,v_{-\lambda_1+\lambda_2}^{\mp}\right]=\pm 3e_{\alpha_2}, &&\left[v_{-2\lambda_1+\lambda_2}^{\pm},\,v_{\lambda_1}^{\mp}\right]=\mp e_{\alpha_3},\\
  &\left[v_{0}^{+},\,v_{0}^{-}\right]=-4h, &&\left[v_{0}^{\pm},\,v_{2\lambda_1-\lambda_2}^{\mp}\right]=\mp 2e_{\alpha_1}, &&\left[v_{0}^{\pm},\,v_{-\lambda_1+\lambda_2}^{\mp}\right]=\mp 2e_{\alpha_3},\\
  &\left[v_{0}^{\pm},\,v_{\lambda_1}^{\mp}\right]=\mp 2e_{\alpha_4}, &&\left[v_{2\lambda_1-\lambda_2}^{\pm},\,v_{-\lambda_1+\lambda_2}^{\mp}\right]=\mp e_{\alpha_4}, &&\left[v_{2\lambda_1-\lambda_2}^{\pm},\,v_{\lambda_1}^{\mp}\right]=\mp 3e_{\alpha_5},\\
  &\left[v_{-\lambda_1+\lambda_2}^{\pm},\,v_{\lambda_1}^{\mp}\right]=\mp 3e_{\alpha_6}.
\end{aligned}\]

\bibliography{refs}

\end{document}